\documentclass{article}

\usepackage[utf8]{inputenc} % allow utf-8 input
\usepackage[T1]{fontenc}    % use 8-bit T1 fonts
\usepackage{hyperref}       % hyperlinks
\usepackage{url}            % simple URL typesetting
\usepackage{booktabs}       % professional-quality tables
\usepackage{amsfonts}       % blackboard math symbols
\usepackage{nicefrac}       % compact symbols for 1/2, etc.
\usepackage{microtype}      % microtypography
\usepackage{parskip}
\usepackage{enumerate}
\usepackage{enumitem}
\usepackage{graphicx}

\title{G\"{o}del's Incompleteness Theorem}
\author{Serafim Batzoglou}

\begin{document}

\maketitle
\begin{abstract}
I present the proof of G\"{o}del's First Incompleteness theorem in an intuitive manner assuming  college-level computational background, while covering all technically challenging steps. I also discuss  G\"{o}del's Second Incompleteness theorem, their connection to G\"{o}del's Completeness theorem, and conclude with brief remarks on implications for mathematics, computation,  theory of mind and AI.
\end{abstract}

\section{Introduction}

I learned about G\"{o}del's Incompleteness theorem as a teenager reading Douglas Hofstadter's "G\"{o}del, Escher, Bach", a unique book that inspired my curiosity for computer science, logic, and philosophy of mind. The theorem has a special place in the history of mathematics.  In a popular list of "top 100" theorems of all time, it ranks 6$^{th}$ after the irrationality of $\sqrt{2}$, Gauss's fundamental theorem of algebra, the countability of the rationals, the Pythagorean theorem and the prime number theorem. And for good reason: to those expecting a foundation of mathematics on a defined collection of axioms from which all theorems are derived, the result must have been as startling as the discovery of $\sqrt{2}$ by the Pythagoreans who believed that all numbers are rational. It has been discussed widely among philosophers and in popular culture for its implications for mathematical foundations, nature of mind and artificial intelligence. 

Despite the theorem's significance, most applied computer scientists are not versed in the theorem's proof or  implications; logic today is much less popular than  statistics and machine learning. In my opinion, the theorem cannot be appreciated unless its proof is understood at a fairly technical level. Informal descriptions, even if well written and precise, leave the reader unclear about the theorem's scope and implications for mathematics and AI.

There are many excellent presentations of the proof (Smith 2013, Smullyan 1991, Nagel and Newman 2001, lecture notes such as by B. Kim). Here, besides having fun writing on the topic, my aim is to present the proof intuitively to computer scientists, assuming minimal college-level background. I omit or push to footnotes many boring details, but cover all technically challenging steps and elements that are pertinent to the theorem's implications. I also discuss briefly G\"{o}del's Second Incompleteness theorem, and the often confusing relation of these theorems to G\"{o}del's Completeness theorem. I conclude with my brief remarks on implications for mathematics, computation, and theory of mind.

\section{Preliminaries}
\subsection{Formal axiomatic systems and Peano Arithmetic}

The theorem concerns formal axiomatic systems for mathematics. Since the time of Euclid, mathematicians have sought to craft formal lists of axioms from which all mathematical truths can be derived. Axiomatic systems usually focus on a specific domain like Euclidean geometry, number theory, or group theory, although some axiomatic systems, such as Zermelo-Fraenkel set theory with the Axiom of Choice, aspire to be rich enough to serve as the foundation for all mathematics. G\"{o}del’s Incompleteness theorem says, in lay language, that any axiomatic system rich enough to capture basic number theory is incomplete. The proof is constructive: it is an algorithm that takes an axiomatic system $\mathcal{L}$ as input, and produces a specific statement $\sigma$ that is either true and unprovable in $\mathcal{L}$, or provable and false, in which case $\mathcal{L}$ is nonsensical.

A formal axiomatic system consists of a \textbf{language} of symbols and syntactic rules with which two types of objects can be composed: \textbf{terms}, which refer to objects in the domain and can be thought of as words of the language,  and \textbf{formulas}, which are mathematical assertions and can be thought of as sentences of the language. Terms and formulas may contain free (input) variables, so their value can vary: terms with no free variables evaluate to a single item in the domain and formulas with no free variables - the \textbf{statements} - evaluate to true or false.\footnote{For example, with $\mathbb{R}$ as the domain  $+, \cdot$ as binary operators and $>, <, =$ as binary predicates, $2.718$, $2\pi r$, and $(x+y)^2$ are examples of terms. $1+1 = 2$, $x > y$, $\forall x \: \forall y \: (x+y)^2 = x^2 + 2xy+ y^2, x < x+1$ are examples of formulas. The first and third are true statements, the second is a formula whose truth depends on $x$ and $y$, and the fourth is a true formula with free variable x.} A finite collection of \textbf{axioms} or axiom schemes (templates that capture a family of axioms) specify all the formulas that can be taken for granted. Then, \textbf{logical rules} specify how to derive a new formula from previous formulas.

For example, Peano Arithmetic (PA) is a simplified language for nonnegative integer arithmetic, with domain $\mathbb{N}$.\footnote{See the section on Completeness - $\mathbb{N}$ is the minimal domain of PA, but larger domains exist.} The symbols for numbers are just $0$ and $S$, and the number $n$ is written as $n$ $S$’s followed by $0$, also denoted as $S^{(n)}0$.\footnote{$S$ formally is a unitary operation whose domain are the terms of the language.} There is a countably infinite set of variable symbols, $x_0$, $x_1$, $x_2$,…. The binary operations $+$, $*$, binary predicates $=$, $<$, $>$ and logical operations $\lor, \land$, $\lnot$, $\forall$, $\exists$ are symbols with their usual meanings. PA's alphabet is the collection of all above symbols,  $\Sigma = \{0, S, +, *, =, <, >, \lnot, \land, \lor, \rightarrow, \forall, \exists, (, ), x_0, x_1, x_2, ... \}$. By convention, we denote variables by $x, y, z,...$; those are not symbols of $\Sigma$.

Terms are defined recursively to be either the basic term $0$, any variable $x$, or $St$, $s+t$, and $s*t$ where $s, t$ are terms. Formulas are defined recursively to be either $s = t, s < t, s > t$ where $s, t$ are terms, or $\lnot \phi, \phi \lor \psi, \phi \land \psi, \forall x \phi, \exists x \phi$, where $x$ is a variable and $\phi, \psi$ are formulas. Parentheses $($, $)$ are used to specify the order in which the above compositions are taken. By convention we denote terms by $s, t, u, ...$ and formulas by $\phi, \chi, \psi, ...$; those are not symbols of $\Sigma$.

The above syntax defines legal terms and formulas of the language. It is simple to distinguish algorithmically between syntactically valid terms, such as $SS0$, $x+(y*Sx)$, syntactically valid formulas, such as $S0 + S0 = SS0$,\footnote{True, as it says $1+1=2$.} $S0 + S0 > SSS0$,\footnote{False, as it says $1+1 > 3$, but syntactically valid.} $\forall y \: (x*Sx > y) \land (x < S0)$, and syntactically invalid strings such as $)(xy\forall S0S$.

PA's axioms come in two groups: (1) domain-specific axioms that specify the properties of natural numbers, addition ($+$), multiplication ($*$), and the relationships $=, <, >$; (2) logical axioms, which represent general logical truths.

Formulas can be derived from axioms or previously derived formulas using three inference rules: (1) from $\phi$ and $\phi \rightarrow \psi$ we can derive $\psi$; this principle is known as \textit{modus ponens} (2) from $\phi$ we can derive $\forall x \phi$ for any variable $x$, whether it appears in $\phi$ or not; this is \textit{universal generalization}\footnote{In some systems it needs to be restricted to variables not free in $\phi$, depending on how it combines with other axioms and inference rules to potentially lead to contradictions.} (3) from $\phi$ we can derive any $\phi'$ that constitutes a renaming of the variables of $\phi$.\footnote{A renaming $\rho$ of variables is such that $\rho(x) \not = \rho(y)$ iff $x$ and $y$ are different variables. This avoids disasters like renaming $x < y \rightarrow \exists z \: (x+z = y)$ to $a < b \rightarrow \exists b \: (a + b = b$).} A \textbf{proof} of a formula $\phi_k$ is simply a list of formulas, $\phi_1, ..., \phi_k$ such that each $\phi_i$ either is an axiom or is derived from previous formulas in the sequence through an inference rule. A formula $\phi$ is a \textbf{theorem} of PA, denoted PA $\vdash \phi$ when there is a proof $[\phi_1, ..., \phi_k = \phi]$ of $\phi$.

The list of PA axioms, logical axioms, and a proof of the statement $S0+S0 = SS0$ are provided in a long footnote.\footnote{There are many equivalent axiomatizations of PA. To simplify things a bit, let us restrict notation to fewer symbols. First, we can do away with the logical symbols $\land, \lor, \leftrightarrow, \exists$ by defining them in terms of $\rightarrow, \lnot, \forall$, as:  $\phi \land \psi := \lnot(\phi \rightarrow \lnot \psi)$; $\phi \lor \psi := \lnot \phi \rightarrow \psi$; $\exists x \phi := \lnot \forall x \lnot \phi$. Then, we can do away with $>, <$: $x > y := \lnot \forall z \: ( x = y+z \rightarrow z = 0)$. A list of axioms that define PA follows (following Weaver, 2014); other possible axiomatizations exist, and studying the axiomatization in detail is not required to understand the main ideas in the incompleteness proof:
\begin{enumerate}
\item Equality: $x = x$;  $\: x = y \rightarrow y = x$; $\: x = y \rightarrow (y = z \rightarrow x = z)$.
\item Successor: $\lnot ( Sx = 0)$; $\: Sx = Sy \leftrightarrow x = y$.
\item Addition: $\: x + 0 = x$; $\: x + Sy = S(x+y)$.
\item Multiplication: $x * 0 = 0$; $\: x * Sy = x+x*y$
\item Induction: $(\phi(0) \land \forall x (\phi(x) \rightarrow \phi(Sx))) \rightarrow \forall x \phi(x)$
\end{enumerate}
Notice that the latter is the familiar induction scheme, and is actually an axiom template from which an infinite collection of axioms derive: for any formula $\phi$, if we show that $\phi(0)$ holds and for every $n$, $\phi(n)$ implies $\phi(n+1)$, then we know that $\phi(n)$ holds for any $n$. Then, the logical axiom schemas comprise:
\begin{enumerate}
\setcounter{enumi}{5}
\item $\phi \rightarrow (\psi \rightarrow \phi)$
\item $(\phi \rightarrow(\psi \rightarrow \chi)) \rightarrow ((\phi \rightarrow \psi) \rightarrow (\phi \rightarrow \chi))$ 
\item $(\lnot \phi \rightarrow \lnot \psi) \rightarrow (\psi \rightarrow \phi)$
\item ($\forall x(\phi \rightarrow \psi)) \rightarrow (\phi \rightarrow \forall x \psi)$, where $\phi$ is a formula that does not contain any free occurrences of $x$.
\item $\forall x (\phi(x)) \rightarrow \phi(t)$ where t is a term that does not share any variables with $\phi$. 
\end{enumerate}
The above axioms are templates that generate infinitely many axiom instances by replacing the formula, term and variable placeholders with any formulas, terms, and variables respectively.
Now we are ready to prove $1+1 = 2$:
(1) $S0+S0 = S(S0+0)$ (axiom 3); (2) $S0 +0 = S0$ (axiom 3); (3) $S(S0 + 0) = SS0 \leftrightarrow S0 + 0 = S0$ (axiom 2); (4) $S(S0+0) = SS0$ (by (2), (3) and modus ponens); (5) $S0+S0 = S(S0+0) \rightarrow (S(S0 +0) = SS0$ $\rightarrow$ $S0+S0 = SS0)$ (axiom 1); (6) $S(S0 +0) = SS0$ $\rightarrow$ $S0+S0 = SS0$ (by (1), (5) and modus ponens); (7) $S0+S0 = SS0$ (by (4), (6) and modus ponens).
}
The key point to understand is that the property of a sequence $\phi_1, ..., \phi_k$ being a proof is syntactic: it can be reduced to simple string manipulation involving the axioms and rules of inference, and as we will see, it can be represented with a PA formula. The collection of theorems of PA is recursively enumerable: it is not hard to write a computer program that mechanically goes through every possible proof and lists all possible theorems of PA.

\subsection{Recursive functions and their representation in Peano Arithmetic}

\subsubsection{Exponentiation in Peano Arithmetic}
This simple language is powerful enough to express number theory. For example, let's write down Fermat's Last Theorem in PA. Recall that the theorem says $x^n + y^n \not = z^n$ for all integers $n>2$ and $x,y,z>0$. However, there is no exponentiation primitive in PA. How do we write $x^n$?

\begin{quote}
\textbf{Notational convention:} PA uses unary representation for numbers. Therefore, $1+3 = 3+1$ is written as $S0 + SSS0 = SSS0 + S0$. In literature, special notation is introduced to distinguish between a number $n$ and its PA representation $\underline{n} = S^{(n)}0 = S....S0$. I think it is OK to follow the convention that a number within a PA formula automatically switches to unary representation. After all, this is also true for decimal: the number $5$ is not the letter $'5'$. A number written on paper automatically switches to a decimal string.
\end{quote}

\noindent Going back to exponentiation, we want a formula $\phi_{EXP}(x, n, y)$ that is true if and only if (\textit{iff}) $y  = x^n$. This is tricky. We want something like $\phi_{EXP}(x, n, y)$ $:=$ $(n = 0 \land y = 1) \lor \exists m \: \exists z \: (n = Sm \land (y = x*z \land \phi_{EXP}(x, m, z))$. Unfortunately, we cannot use $\phi_{EXP}$ in its own definition.

The trick is to pack all of the intermediate exponents $x, x^2, ..., x^n$  into a single big number $N$. Then we can devise a formula $\phi_{EXP}$ that "talks" about $N$ having the property of packing all the intermediate exponents of $x$. In principle it is no surprise that we can pack an arbitrary list of numbers into a single big number. In practice it is tricky to do so with the limited syntax of PA. 

The ability to pack a list of numbers $k_1, ... , k_n$ into a single number is a "subroutine" of the incompleteness proof, and is the key insight allowing PA to represent recursive functions, as described in the next section. The original way G\"{o}del accomplished this was using the ancient Chinese Remainder Theorem, and is still one of the cleanest ways:

\begin{quote}
\textbf{The Chinese Remainder Theorem:} Given $k_1, \ldots k_n \in \mathbb{N}$ and pairwise coprime $m_1, \ldots , m_n \in \mathbb{N}$ with $m_i > k_i, 1 \leq i \leq n$, there is a unique $N < m_1 m_2 \ldots m_n$ such that for all $1\leq i \leq n$, $k_i \equiv N$ mod $m_i$. 
\end{quote}

\noindent The theorem seems made-to-order for packing lists into a single integer. Given any $k_1, ..., k_n$, we can pack them into a single $N$ by fixing a list of coprime $m_i > k_i$, finding the unique $N < m_1 m_2 \ldots m_n$ that "works", and then retrieving any $k_i$ by taking $N$ mod $m_i$. One important subtlety is that the $m_i$ have to be describable in a compact way - say, with at most a single additional integer on top of $N$. Let's say we can do that. Then, intuitively, there are $m_1 m_2 \ldots m_n$ possible remainder lists $[k_1, ..., k_n]$ of dividing a number $N$ by each of the $m_i$ in order. Because the $m_i$ are pairwise coprime, by letting $N$ cycle through all values $N = 1, ..., m_1 m_2 \ldots m_n$, all possible such remainder lists are visited.\footnote{The cyclic group $\{1, \ldots, m_1m_2 \ldots m_n\}$ under multiplication is isomorphic to the product of cyclic groups $\{1, \ldots m_1\} \times \cdots \times \{1, \ldots, m_n\}$ as long as the $m_i$ are pairwise coprime.}

To implement this in PA, first we define the $mod$ function $c \equiv a$ mod $b$: $\phi_{MOD}(a, b, c) := \exists n$ $(a = b*n + c \land c <b)$.

Next, we need a way to pick specific $m_i > k_i$ to implement the packing. The $m_i$ have to be pairwise coprime and collectively describable with a single number. We can accomplish that by picking a suitable $b > \max\{k_1, ..., k_n\}$, and letting $m_i = 1 + ib$. As long as every number among $2, ..., n-1$ divides $b$, the $m_i$ are guaranteed to be pairwise coprime.\footnote{\textbf{Proof:} Let $p \mid 1+ib$, $p\mid 1+jb$, $1\leq i < j \leq n$. Then $p \mid (j-i) b$ so $p \mid j-i$ or $p \mid b$. Also, $j-i \mid b$ by assumption, therefore $p \mid b$. Then $p \mid 1+ ib$ implies $p \mid 1$ so $p = 1$.}

Let $b$ be a multiple of $n!$ that is greater than $\max\{k_1, ..., k_n\}$. Let $m_i = 1+ib$. Let $N$ be the unique integer $< m_1 m_2 \ldots m_n$ such that $k_i$ is the remainder of $N$ over $1 + ib$ for all $1 \leq i \leq n$. $N$ and $b$ are computable given any $k_1, ..., k_n$.

Then, given the packed number $N$ and also $b, i$, we can extract $k_i$ as the remainder of $N / (1+ib)$.

Finally we are ready to define exponentiation by packing $x, x^2, ..., x^{n-1}$ into $N$. In the "code" below notice the use of the "subroutine" $\phi_{POW-N}$, to ensure that $N$ is minimal.\footnote{This is  not strictly needed: the definition would still be valid without ensuring that $N$ is minimal.}
$$\phi_{POW-N}(N, x, n, y) := \{n = 0 \land y = S0\} \lor  \{\exists m \: (n = Sm) \land$$
$$[\exists b \: \forall i<n \: \exists j \: (b = i * j) \land\phi_{MOD}(N, Sb, x) \land \phi_{MOD}(N, S(n*b), y) \land$$ 
$$(\forall i<m \: \forall z \:  \phi_{MOD}(N, S(Si*b), z) \rightarrow \phi_{MOD}(N, S(SSi * b), x*z))]\}$$
$$\phi_{POW}(x, n, y) := \exists N \: \phi_{POW-N}(N, x, n, y) \land \forall M \: (\phi_{POW-N}(M, x, n, y) \rightarrow N \leq M)$$

\noindent In the above formula, $\{\}, []$ are used instead of $()$ just for clarity; these are not symbols of PA. Also, "$\forall x (x<y) \rightarrow \phi$" is abbreviated to $\forall x<y \: \phi$. The formula is a mouthful, but once we have it we can use it as a subroutine. 

Fermat's Last Theorem can be expressed as follows:
$$\forall n\: (n>SS0) \rightarrow [\forall x\: \! >\! 0 \: \forall y\: \! >\! 0 \: \forall z\: \! >\! 0 \: \lnot \exists a \: \exists b \: \exists c \: (a+b =c \land \phi_{POW}(x, n, a) \land \phi_{POW}(y, n, b) \land \phi_{POW}(z, n, c))]$$
\subsubsection{Recursive Functions}

Exponentiation is an example of a recursive function, the class of functions that can be calculated by an effective method, which are precisely the functions computable by Turing machines according to the Church-Turing thesis. Recursive functions defined on every input on $\mathbb{N}$ are the functions computable by a Turing machine that is guaranteed to halt. PA is rich enough to \textbf{represent} these functions. In particular:

\begin{quote}
$f: \mathbb{N} \rightarrow \mathbb{N}$ is \textbf{represented} by formula $\phi(x, y)$ iff for all $x\in \mathbb{N}$, PA $\vdash \forall y \: \phi(x,y) \leftrightarrow y = f(x)$.
\end{quote}

\noindent Representing a function $f$ is a strong statement: not only we want a formula $\phi(x,y)$ that is true iff $y = f(x)$, but additionally, we want the statement $\forall y \: \phi(x,y) \leftrightarrow y = f(x)$ to be derivable in PA.

A couple of things to note: first, in the above definition, $f(x)$ is a numerical value represented by $S^{(f(x))}0$; there is no symbol for $f$ in PA and a separate such statement is derived for every $x$. Second, the  definition generalizes to functions $f: \mathbb{N}^k \rightarrow \mathbb{N}$ represented by formulas $\phi(x_1,..., x_k, y)$.

Recursive functions are defined inductively on $\mathbb{N}^k \rightarrow \mathbb{N}$ starting from a number of basic functions that can be combined through a few composition rules. There are many alternative, equivalent definitions, and excellent references are available online (Dean W, 2020).

\begin{quote}
\textbf{Basic recursive functions.} The basic recursive functions are:
\begin{itemize}
\item The \textit{zero} function on $\mathbb{N}^k \rightarrow \mathbb{N}$, $0_k(x_1, ..., x_k) = 0$.
\item The \textit{successor} function on $\mathbb{N} \rightarrow \mathbb{N}$, $\mathcal{s}(x)  = x+1$. 
\item The \textit{projection} functions on $\mathbb{N}^k \rightarrow \mathbb{N}$, $1 \leq i \leq k$, $\pi_i(x_1, ..., x_k) = x_i$.
\end{itemize}
A function $f: \mathbb{N}^k \rightarrow \mathbb{N}$ is \textbf{primitive recursive} if it is obtained from the primitive functions by finitely many applications of the following two compositional rules:
\begin{itemize}
\item \textit{Composition}: $f(x_1, ..., x_k) = g(h_1(x_1, ..., x_k), ..., h_l(x_1, ..., x_k))$, 
\subitem where $g: \mathbb{N}^l \rightarrow \mathbb{N}$, $h_i: \mathbb{N}^k \rightarrow \mathbb{N}$ are functions constructed at a previous stage.
\item \textit{Primitive Recursion}:
\subitem  $f(x_1, ..., x_{k-1}, 0) = h(x_1, ..., x_{k-1})$ 
\subitem  $f(x_1, ..., x_{k-1}, n+1) = g(x_1, ..., x_{k-1}, n, f(x_1, ..., x_{k-1}, n))$, 
\subitem where $g: \mathbb{N}^{k+1} \rightarrow \mathbb{N}$, $h: \mathbb{N}^{k-1} \rightarrow \mathbb{N}$ are functions constructed at a previous stage.
\end{itemize}

\end{quote}

\noindent A rich collection of functions can be constructed with the above primitives. For example, constant functions, addition, multiplication and exponentiation can be defined through primitive recursion:
\begin{quote}
\begin{itemize}
\item \textit{Constant} $K$: $K(x) = \mathcal{S}...\mathcal{S}(0_1(x)...)$ ($K$ times)
\item \textit{Addition:} $+(x, 0) = x$; $+(x, y+1) = \mathcal{S}(+(x, y))$ 
\item \textit{Multiplication:} $*(x, 0) = x$; $*(x, y+1) = +(x, *(x, y))$ 
\item \textit{Power:} $x^0 = Pow(x, 0) = 1$; $x^{y+1} = Pow(x, y+1) = *(x, Pow(x, y))$ 
\end{itemize}
\end{quote}

\noindent However, this does not exhaust all computable functions. We need one additional operation, that of minimization. Consider a computable function $g(x,y)$ such that for every $x$ there is some $y$ such that $g(x, y) = 0$. Then, we can find the minimum such $y$ by computing $g(x, y)$ for all values $y = 0, 1, 2, ...$ until we hit $g(x, y) = 0$. Therefore we want to say that minimizing for $y$ this way is also computable:
\begin{quote}
\begin{itemize}
\item \textit{Minimization:} $f(x) = \mu_{y}g(x,y) := min_y$ s.t. $g(x,y) = 0$
\subitem where $g: \mathbb{N}^2 \rightarrow \mathbb{N}$ is recursive and $\forall x \: \exists y \: g(x, y) = 0$
\end{itemize}
\end{quote}

\noindent A function $f$ is \textbf{recursive} if it is obtained from the primitive functions through finite applications of composition, primitive recursion and minimization. If $f$ is defined on all inputs (guaranteed if the minimization applications are defined on all inputs) then it is \textbf{total recursive}. 

The recursive functions are precisely the functions computable by a Turing machine. The total recursive functions are precisely the functions computable by a Turing machine that halts on all inputs. There exist total recursive functions that are not primitive recursive such as the Ackermann function, but they are hard to come up with.

\subsubsection{Representing recursive functions in PA}

Peano Arithmetic is rich enough to represent recursive functions. Recall that a function $f(x)$ is represented by $\phi(x,y)$ if for all $x$, $\vdash \forall y \: \phi(x,y) \leftrightarrow y = f(x)$. The notion readily generalizes to functions of more than one argument. To show that PA can represent recursive functions, we can go case-by-case in the definition above. Most cases are easy to show, and formal proofs will be omitted here:
\begin{quote}
\begin{itemize}
\item Zero: $\phi_0(x_1, ..., x_k, y) \: := \: y = 0$.
\item Successor: $\phi_S(x, y) \: := \: y = S(x)$.
\item Projection: $\phi_{\pi_i}(x_1, ..., x_k, y) \: := \: y = x_i$.
\item Composition:  assume $\psi_i(x_1, ..., x_k, y)$ represents $h_i$ for $1\leq i \leq l$ and $\chi(x_1, ..., x_k, y)$ represents $g$, then 
 $$\phi_{g \circ [h_1, ..., h_l]}(x_1, ..., x_k, y) \: := \: \exists y_1 ... \exists y_k (\psi_1(x_1, ..., x_k, y_1) \land ... \land \psi_k(x_1, ..., x_k, y_k) \land \chi(y_1, ..., y_k, y))$$ represents $f$.
\item Minimization: assume $\psi(x, y, z)$ represents $g$, then $\phi_{\mu}(x, y) \: := \: \psi(x, y, 0) \land \forall w (w < y \rightarrow \lnot \psi(x, w, 0))$ represents $f(x) = \mu_y g(x, y)$.
\end{itemize}
\end{quote}

\noindent The difficult case is primitive recursion. Consider a function $f$ defined in terms of recursive functions $g, h$ by $f(x, 0) = h(x)$ and $f(x, n+1) = g(x, n, f(x, n))$.\footnote{This easily generalizes to multiple variables; see also: \\http://www.michaelbeeson.com/teaching/StanfordLogic/Lecture11Slides.pdf.} How do we represent it with a formula?

What we would want is: $y = f(x, n) \leftrightarrow \exists y_0 ...\exists y_n \: y_0 = h(x) \land y_1 = g(x, 0, y_0) \land ... \land y_n = g(x, n-1, y_{n-1}) \land y = y_n$. But this is not a valid formula of PA.

So, what we really want is to be able to refer within a formula to an entire sequence of values, $y_0, ..., y_n$ in a way that relates them to one another. We have already encountered this situation in defining $\phi_{POW}$, which is representing the recursive function $x^y = x*x^{y-1}$, whereby we needed to represent $1, x, x^2, ..., x^y$ within a single number. The same trick works in the general case of primitive recursion.

We first define a function  $\beta(N, i)$ representable in PA, which encodes arbitrary sequences of numbers $n_1, ..., n_k$ into a single number $N$, whereby $\beta(N, i) = n_i$ for $1 \leq i \leq k$.

\begin{quote}
\textbf{The $\beta$-function lemma. } There is a function $\beta: \mathbb{N}^2 \rightarrow \mathbb{N}$ representable in PA s.t. for any $N,i$, $\beta(N,i) < N$, and for any sequence $n_1, ..., n_k \in \mathbb{N}$ there exists $N\in \mathbb{N}$ such that $\beta(N, i) = n_i$ for all $1\leq i \leq k$. 
\end{quote}

\noindent This function appears in the original proof by G\"{o}del. It is implemented through the Chinese Remainder Theorem similarly to $\phi_{POW}$ above; I am pushing details to a footnote.\footnote{The $\beta$ function is defined so as to ensure that for any $n_1, ..., n_k$ there is an $N$ such that $\beta(N, i) = n_i$ for $1\leq i \leq k$. Without loss of generality let $n_0 = k$ be the length of the list.

Let $c = \max\{n_0, ..., n_k\}$. let $b=c!$ and observe that $b+1, ..., (k+1)b + 1$ are relatively prime. By the Chinese Remainder Theorem, there is a unique $n < \Pi_{i = 1}^{k+1}  (ib+1)$ such that $n \equiv n_i$ mod $(ib + 1)$ for all $1\leq i \leq k+1$.

Now from $n, b$ we can retrieve $k = n$ mod $(b+1)$ and then retrieve all remaining $n_1, ..., n_k$. To further compress to a function $\beta$ of two arguments as specified by the lemma, define a pairing function $\pi: \mathbb{N}^2 \rightarrow \mathbb{N}$ and two projection functions $\pi_L, \pi_R: \mathbb{N} \rightarrow \mathbb{N}$ such that for any $i, j \in \mathbb{N}$ we have $\pi_L(\pi(i, j)) = i$ and $\pi_R(\pi(i,j)) = j$, with $\pi, \pi_L, \pi_R$ primitive recursive. The specific implementation is not important, but $\pi(i, j) = (i+j)^2 + i + 1$ works because it is injective and easy to represent in PA.

Define $N := \pi(n, b)$ and $\beta(N, i) := \pi_L(N)$ mod $((i+1)\pi_R(N) + 1) = n_i$ for $0 \leq i \leq k$. To define primitive recursion, the assertion will be that $\exists N$ such that $\beta(N, i)$ satisfy the recursive definition for $k, n_1,..., n_k$, and moreover that such $N$ is the minimum satisfying the definition.}

All the above functions can be represented in PA using the previously defined recursion formulas (zero, successor, projection, composition and minimization). The details are omitted (see for instance lecture notes by B. Kim).
\begin{quote}
\begin{itemize}
\item Primitive Recursion: let $\phi_g, \phi_h$ represent $g$ and $h$, and $\phi_{\beta}(N, i, y)$ represent $\beta(N,i)$; then $f(x,n)$ is represented with $\phi_{REC}(x, n, y)$: 
\end{itemize}
\end{quote}
$$\phi_{REC-N}(N, x,n,y) \: :=  (\exists y_0 \: \phi_{\beta}(N, 0, y_0) \land \phi_h(x, y_0)) \land \forall i \: [i<n \rightarrow $$ $$\exists y_i \: \exists y_{i+1} \: \phi_{\beta}(N, i, y_i) \land \phi_{\beta}(N, Si, y_{i+1}) \land \phi_g(x, i, y_i) \land \phi_g(x, Si, y_{i+1})] \land \phi_{\beta}(N, n, y)$$
$$\phi_{REC}(x, n, y) := \exists N \: \phi_{REC-N}(N, x, n, y) \land \forall M \: (\phi_{REC-N}(M, x, n, y) \rightarrow N \leq M)$$

\noindent Notice the "subroutine" $\phi_{REC-N}$ within $\phi_{REC}$ so as to ensure that  $N$ is minimal - this condition is actually not needed for the definition to work, but it's cleaner this way. Once the above formulas are laid out, it is tedious to demonstrate that in each case it can be derived within PA that $y = f(...)$ iff $\phi(..., y)$, i.e., to prove PA $\vdash \forall y \: (\phi(..., y) \leftrightarrow y = f(...))$. This ability of PA is at the heart of the incompleteness theorem. This makes PA as expressive as recursive functions or general Turing machines. Then, the incompleteness theorem follows just as surely as the undecidability of the halting problem does for Turing machines. It is instructive to ask what is a (\textit{the}) key property of PA that makes it so. In fact, it is the combination of  addition \textit{and} multiplication. Presburger arithmetic is essentially PA without multiplication, and Skolem arithmetic is PA without addition; both are complete and decidable.\footnote{There is a doubly exponential procedure for Presburger arithmetic and a triply exponential procedure for Skolem arithmetic, which determine whether a formula is a theorem. It is notable that even slightly more power turns the theories undecidable. For example, extending Skolem arithmetic with the successor operator $S$ enables implementing addition; extending with the $<$ predicate allows implementing the successor operator, $S(x) := y>x \land \forall z>x \: (z = y \lor z>y)$ and therefore also renders the resulting theory undecidable.} One can see why: key to the ability to pack an arbitrary list of numbers $[n_1, ..., n_k]$ into a single number is the combination of addition and multiplication. In the particular implementation of the $\beta$ function as described here, the Chinese Remainder Theorem is used. However, this is just an implementation choice: other choices are possible, such as picking the largest prime number $p > \max_i\{n_1, ..., n_k\}$ and packing all numbers into the pair $(n,p)$ where $n = \Sigma_{i=1}^k n_i p^{i-1}$. Regardless, without both addition and multiplication, this can't take place.\footnote{Weaker systems than PA that exhibit incompleteness are possible. Primitive Recursive Arithmetic (PRA) is a well-known weaker system that differs from PA in the following ways: (1) there are no quantifiers; (2) instead, there is a separate symbol for each primitive recursive function, which is now needed because of lack of quantifiers; (3) there are only two inference rules, modus ponens and variable substitution - the quantifier rule is dropped. Successor, addition and multiplication are preserved. PRA, and even slightly weaker subsystems allow the incompleteness proof to go through. Another example is the theory of hereditarily finite sets (Swierkowsky 2003).}

\section{Theorem Statement}

The technical statement of the theorem and its slight variants can be a bit cumbersome. Informally, the theorem states:

\begin{quote}
Any reasonable extension $\mathcal{L} $ of PA is either inconsistent or incomplete. Specifically, there exist a statement $\sigma$ such that either $\sigma$ is $true$ and $\mathcal{L}  \not \vdash \sigma$, or $\mathcal{L}  \vdash \bot$.
\end{quote}

\noindent The above statement is blatantly informal. First, what is a "reasonable" extension of PA? That would be a language $\mathcal{L} $ that includes all of PA, plus at most countably many additional symbols, and a recursive set of additional axioms or axiom templates, so that checking the validity of a proof in $\mathcal{L} $ is a primitive recursive task. Any theorem of PA would be a theorem of $\mathcal{L} $.

A few remarks:

\begin{itemize}
\item In the original proof, an additional requirement was that the axiomatic system is $\omega$-\textit{consistent}: that there is no formula $\phi(n)$ such that for any $n$, $\mathcal{L}  \vdash \phi(n)$, while simultaneously $\mathcal{L}  \vdash \exists n \lnot \phi(n)$. However, $\omega$-consistency is too strong an assumption. Soon after G\"{o}del's original proof, Rosser strengthened it by introducing a trick explained below that removes this assumption. This newer version is what is also known as the G\"{o}del-Rosser theorem.
\item This is the first of two famous incompleteness theorems. The Second Incompleteness theorem informally states that $\mathcal{L} $ cannot prove its own consistency. 
\item Technically different versions of the theorem can be proven through computability theory, relating to the Halting problem, or through Kolmogorov complexity, relating to the smallest possible program that outputs a given string. I will discuss these angles briefly below.
\end{itemize}

\section{Proof Outline}
\subsection{Main idea - arithmetization and diagonalization}
The main idea of the proof is based on the liar's paradox "this sentence is false", modified within PA to say "this statement is not provable". Self-referential statements are well known troublemakers in foundational mathematics, logic and philosophy. The most famous example is perhaps Russel's paradox in naïve set theory of the set $x$ of all sets that are not members of themselves. Is $x$ a member of itself? It can neither be, nor not be. Or Russel's barber "\textit{who shaves all those, and those only, who do not shave themselves}".  Does the barber shave himself? 

Similarly here, the statement $\sigma :=$ "this statement is not provable" is either true and not provable, or provable and false. That implies that PA is either incomplete or inconsistent. The proof is an algorithm: starting from the symbols and axioms of PA, it generates a statement $\sigma$ such that "$\sigma$ is true if and only if $\sigma$ is not provable" is a theorem of PA. The algorithm is general enough to be applicable to any reasonable extension of PA. How can such a $\sigma$ be constructed? Through two key steps in the proof: \textbf{arithmetization} and \textbf{diagonalization}.

\textbf{Arithmetization.} If we could devise a PA formula $\phi_{PROVABLE}$ that applies to other formulas, and a statement $\sigma$ such that $\sigma \leftrightarrow \lnot \phi_{PROVABLE}(\sigma)$ is a theorem of PA, we would be done. Unfortunately, PA's domain is numbers and not formulas. And here comes the proof's first major ingenuity: define a function that \textit{arithmetizes} formulas and proofs of formulas. Seen as strings of letters, formulas and proofs can be encoded into numbers in a variety of ways such as ASCII. Today this may seem obvious; in 1931 when the theorem was proven, it was brilliant.

Having encoded formulas into numbers, $\sigma$ is mapped to a number, call it $\lceil \sigma \rceil \in \mathbb{N}$. Lists of formulas are also mapped to numbers, $\lceil [\phi_1, ..., \phi_k] \rceil \in \mathbb{N}$. A proof of $\sigma$ in PA is simply a list $[\phi_1, ..., \phi_k = \sigma]$ such that each $\phi_i$ is either an axiom or a consequence of preceding formulas under an inference rule. The latter is a syntactic property: it can be performed on the proof text by a well defined procedure in finite steps. The same procedure translates to arithmetic operations on a number $m$ representing the proof $[\phi_1, ..., \phi_k = \sigma]$. This can represented in PA, as we will see below, with $\phi_{PROOF-OF}(m, \lceil \sigma \rceil)$.\footnote{Note that the arguments of all these formulas are numbers and are actually represented in unary $S...S0$ within PA.} Then, the formula $\phi_{PROVABLE}(\lceil \sigma \rceil)$ is expressed as $\exists m \: \phi_{PROOF-OF}(m, \lceil \sigma \rceil)$.\footnote{This latter property is \textit{not} decidable; it is the only non-decidable property defined in PA in the incompleteness proof.}  If we can now devise a $\sigma$ such that $\sigma \leftrightarrow \lnot \phi_{PROVABLE}(\lceil \sigma \rceil)$, we are done. This is accomplished through \textit{diagonalization}. 

\textbf{Diagonalization.} G\"{o}del described two different diagonalization methods in his original paper. An intuitive, informal one, and a  rigorous one in full formal detail. The intuitive method goes as follows. Consider all formulas of PA with a single free variable (plus optionally other bound variables), and moreover require the free variable to be specifically $x_0$. These formulas can be ordered lexicographically $\phi_1(x_0)$, $\phi_2(x_0)$, $\phi_3(x_0)$,... in a manner that is computable thus representable in PA. Therefore we can talk in PA about the $n^{th}$ such formula, $\phi_n(x_0)$, as long as $n$ itself is definable. For any definable $n,m \in \mathbb{N}$, the statement $\phi_n(m)$ is definable in PA, contains no free variables and is either true or false. Consider now the formula $\phi_{NOT-PROVABLE-SELF}(x_0) :=$ "the statement $\phi_{x_0}(x_0)$ is not provable".\footnote{To construct this step-wise, first define $\phi_{NOT-PROVABLE}(x_0, x_1) :=$"the statement $\phi_{x_1}(x_0)$ is not provable". This is a formula of two variables. Letting $x_0 = m, x_1=n$, it says that the $n^{th}$ formula in the lexicographic list of single-variable formulas, when given input $m$ becomes a statement that is not provable. Then define $\phi_{NOT-PROVABLE-SELF}(x_0) := \phi_{NOT-PROVABLE}(x_0,x_0)$.} This is definable in PA and states that the $x_0^{th}$ formula in the list, with number $x_0$ as argument is not provable. The formula $\phi_{NOT-PROVABLE-SELF}(x_0)$ itself has the single free variable $x_0$, and therefore is on the list say at position $M$.\footnote{$M$ is definable in PA as it can be calculated from the syntax of $\phi_{NOT-PROVABLE-SELF}$.} Now consider what is the value of $\phi_{NOT-PROVABLE-SELF}(M)$. This states "the statement $\phi_{NOT-PROVABLE-SELF}(M)$ is not provable". If it is true, then the statement is not provable, therefore PA is incomplete. If it is false, then the statement is provable and false, so PA is inconsistent. The argument is called "diagonalization" because we apply the $M^{th}$ formula to the argument $M$. In the proof that follows, a similar argument is used except with a different arithmetization, similar to G\"{o}del's rigorous arithmetization.

\subsection{Outline of main steps}
The proof takes the following main steps:
\begin{enumerate}
\item \textbf{Encode formulas into numbers.} First, a function $G$ is defined that maps strings into numbers in $\mathbb{N}$. Any formula $\phi$ maps to $G(\phi)$, denoted as $\lceil \phi \rceil \in \mathbb{N}$, from which $\phi$ can be retrieved. This way, formulas can "talk" about other formulas. For example, we can write a formula $\phi_{FORMULA}(n) =$ "$n$ is the encoding of a syntactically valid formula". 
\item \textbf{Extend the map to lists of formulas.} The mapping $G$ extends to ordered lists of formulas, $G([\phi_1, ..., \phi_k]) = \lceil [\phi_1, ..., \phi_k] \rceil \in \mathbb{N}$. Given a number $n = \lceil [\phi_1, ..., \phi_k] \rceil$, the list of formulas can be retrieved.
\item \textbf{Express the notion of a provable formula.} A proof is just an ordered list of formulas with a special property: every formula in the list is either an axiom or a consequence of previous formulas through one of the three derivation rules. This is a syntactic property, which is easily seen to be computable and is now captured in a formula $\phi_{PROOF-OF}(m, n) =$ "$m$ encodes a proof of the formula encoded by $n$", which is true and derivable in PA iff $m = \lceil [\phi_1, ..., \phi_k] \rceil $, $n=\lceil \phi_k \rceil$, and $[\phi_1, ..., \phi_k]$ is a proof of $\phi_k$. Then the notion of a provable formula can be captured in a formula: $\phi_{PROVABLE}(n) =$ "$n$ encodes a provable formula" $=\exists m$ $\phi_{PROOF-OF}(m, n)$. 
\item \textbf{Devise a self-referential formula by a variable-substitution trick.} Then, a rather complex formula $\sigma$ is expressed, and it is shown within PA that $\sigma \leftrightarrow \lnot \phi_{PROVABLE}(\lceil \sigma \rceil)$. This is known as the "fixed point lemma" and is a tricky and fun part of the proof. First, the "diagnonalization" $\mathcal{D}$ of a formula is defined as a procedure that takes a number $m$, and if $m = \lceil \phi \rceil$, substitutes the zero or more occurrences of the specific variable $x_0$ within $\phi$ with the number $m$ (expressed $S^{(m)}0$): $\phi(x_0, y, z, ...)$ turns into $\phi(m, y, z, ...)$. Then, $\mathcal{D}$ returns $\lceil \phi(m, y, z, ...) \rceil$. This recursive function is represented in PA with a formula $\psi_{DIAG}(x, y)$, and using this formula, statement $\sigma$ is constructed (the punchline will be reserved for the respective section below).
\item \textbf{Conclude that the resulting formula is either true and not provable, or false and provable.} $\sigma$ is a numerical statement of no free variables, and is either true or false in $\mathbb{N}$. If $\sigma$ is provable in PA, then PA $\vdash \lnot \phi_{PROVABLE}(\lceil \sigma \rceil)$, a contradiction.\footnote{A proof of $\sigma$ in PA can be encoded with into a number $m^*$ to derive $\phi_{PROOF-OF}(m^*, \lceil \sigma \rceil)$ in PA, which implies $\phi_{PROVABLE}(\lceil \sigma \rceil)$.} If on the other hand $\lnot \sigma$ is provable, then PA $\vdash \phi_{PROVABLE}(\lceil \sigma \rceil)$. Then, assuming that PA is $\omega$-consistent (i.e., it cannot derive $\psi(n)$ for every $n \in \mathbb{N}$ and simultaneously derive $\exists n \lnot \psi(n)$), this is a contradiction. This latter assumption of $\omega$-consistency is a subtle point: it can be removed by slightly modifying the above proof as will be discussed below, with a simple trick devised by J. Barkley Rosser after the original proof appeared. The conclusion is that $\sigma$ is true and not provable in PA, therefore PA is incomplete. 
\end{enumerate}
\medskip

\noindent \textbf{What is needed for the proof to go through.} Let us list a few key elements needed for the above proof to go through for an axiomatic system $\mathcal{L}$, which in our case is PA:
\begin{enumerate}
\item $\mathcal{L}$'s domain has to be rich enough for $\mathcal{L}$'s formulas to be encoded unambiguously to the domain, so that $\mathcal{L}$ can "talk" about its own formulas. 
\item $\mathcal{L}$ has to be expressive enough to express the notion of "provable formula" with a formula, $\phi_{PROVABLE}$ and also the diagonalization procedure that involves variable substitution. Generally, what is needed is a system that can represent all computable functions. 
\item $\mathcal{L}$'s axioms have to be powerful enough for steps 4 and 5 to go through. A proof of $\sigma$ in $\mathcal{L}$, call it $[\phi_1, ..., \phi_k = \sigma]$ will only lead to a contradiction with the assertion $\lnot \phi_{PROVABLE}(\lceil \sigma \rceil)$ if $\phi_{PROOF-OF}(\lceil [\phi_1, ..., \phi_k]\rceil, \lceil \phi_k \rceil)$ is guaranteed to be provable in $\mathcal{L}$ whenever $[\phi_1, ..., \phi_k]$ is a valid proof. That is to say, $\mathcal{L}$'s axioms should be powerful enough to check the validity of a fully spelled-out proof. Again, what we need here is the power to represent computable functions. 
\end{enumerate}
\medskip

Note that the notion "provable" is not decidable. One cannot construct a program that takes a formula as input, churns away, and always returns a correct "yay" or "nay" as a result. Such a program would solve the halting problem. However, \textit{given} a proof of a formula in sufficient detail, it is computable to check that the proof is correct. Therefore, a proof can be encoded by $G$ into a number $m$ and utilized to derive $\phi_{PROOF-OF}(m, \lceil \sigma \rceil)$, which immediately implies $\phi_{PROVABLE}(\lceil \sigma \rceil)$. 

These key elements can be satisfied in multiple ways. The language does not have to be PA. There is great freedom in designing the arithmetization. The self-referential formula can be changed substantially. The proof can be thought of as a pseudocode implementation of a  more general idea.

\subsection{More details on main steps}

\subsubsection{Step 1: Arithmetize formulas}
The formulas of PA are a countable, primitive recursive collection of strings of letters that can be mapped 1-1 to $\mathbb{N}$. There are many possible mappings, such as for example ASCII.

G\"{o}del provided a clever mapping before the time of ASCII, based on exponentiating prime numbers to powers that codified letters in the alphabet of PA. First, codify each symbol of the alphabet: $\{0:1, S:2, +:3, *:4, =:5, \lnot:6, \land:7, \lor:8, \rightarrow:9, \forall:10, (:11, ):12, x_0 :13, x_1 :17, x_2 :19, ... \}$. Then, convert a string $s = s_1 s_2... s_k$ of length $k$ to a list of numbers $[n_1, ..., n_k]$ using the above code. Then pick the first $k$ prime numbers in $\mathbb{N}$, $2, 3, 5, ..., p_k$, and calculate the number $2^{n_1} 3^{n_2} \ldots p_k^{n_k}$ that uniquely represents $s$, from which $s$ can be retrieved using prime factorization. For example, $S0 + S0 = SS0$ is converted to $[3,1,5,3,1,9,3,3,1]$ and then to $2^3 3^1 5^5 7^3  11^1  13^9  17^3  19^3  23^1$. It is a big number, but since nobody will actually ever implement this algorithm, big numbers are cheap.

Given any formula $\phi$, we can thus compute a unique $\lceil \phi \rceil \in \mathbb{N}$. Conversely, given any number $n$, we can compute whether $n$ is the G\"{o}del number $\lceil \phi \rceil$ of a syntactically valid formula $\phi$.

\subsubsection{Step 2: Arithmetize lists of formulas}

Sequences of formulas can  be encoded by using the $\beta$ function: $\lceil [\phi_1, ..., \phi_k]\rceil$ is defined for any list of numbers to be the minimum $N$ such that $\beta(N, 0) = k$ and $\forall_{1\leq i \leq k} \: \beta(N, i) = \lceil \phi_i \rceil$.\footnote{Terms and formulas can also be arithmetized thus according to their structure. For example, $s=t$ becomes $\lceil [\lceil = \rceil, \lceil s \rceil, \lceil t \rceil]\rceil$, $\forall x \phi$ becomes $\lceil [ \lceil \forall \rceil, \lceil x \rceil, \lceil \phi \rceil ] \rceil$ and so on. This encoding makes it much easier to write "code" in this language. See lecture notes by B. Kim.}

A proof of a formula $\phi$ is a finite list of formulas, $[\phi_1, …, \phi_k = \phi]$, such that every $\phi_i$ for $1 \leq i \leq k$ is either an axiom or is derived from $\phi_1, …, \phi_{i-1}$ using PA's rules of inference. Given $n \in \mathbb{N}$, we can compute whether $n$ encodes a syntactically valid sequence of formulas, and whether the sequence constitutes a valid proof in PA: we just need to check that every $\phi_i$ is either an axiom\footnote{This can be checked syntactically through a finite number of substitutions of variables within the axioms with terms appearing in $\phi_i$.} or a consequence of some preceding $\phi_j, \phi_l, j,l<i$ in the sequence through one of the three rules of inference.

\subsubsection{Step 3: Express the notion of a provable formula}
Now that formulas, sequences of formulas and proofs are mapped into numbers, PA can deal with them. Every recursive function is representable in PA, therefore any computable property $P(x)$ where $x$ is a formula or a proof can be first mapped to the indicator function $f_P(\lceil x \rceil)$ and then represented in PA with a formula $\phi_f(\lceil x \rceil, y)$ such that the statement $y = f_P(\lceil x \rceil) \leftrightarrow \phi_f(\lceil x \rceil, y)$ is a theorem of PA. Since $y$ takes values $0$ or $1$ in this case, we might as well let $\phi_f' (\lceil x \rceil) := \phi_f(\lceil x \rceil, S0)$, which is true iff $P(x)$. 

Let's see some examples:  

The property $P(x)$ could be "$x$ is a syntactically valid formula in PA". This is certainly a computable property. We can turn this into a PA formula that checks syntactic validity of encoded formulas:

\begin{quote}
$\phi_{FORMULA} (n) := \:$ "$n = \lceil \psi \rceil$ for some syntactically valid formula $\psi$ of PA"
\end{quote}

\noindent Also, checking whether a formula has a single, specific free variable $x$ is computable:
\begin{quote}
$\phi_{SVF-x} (n) := \:$ "$\phi_{FORMULA} (n)$, $n = \lceil \psi \rceil$, and $\psi = \psi(x)$ has the single free variable $x$"\footnote{We can define such a formula for any specific $x$, or a formula checking for any single variable, or for two variables, etc. All these syntactic properties are easily shown to be computable and therefore representable in PA.}
\end{quote}

\noindent Checking whether a sequence of formulas is a valid, complete proof within PA of the last formula of the sequence is computable:
\begin{quote}
$\phi_{PROOF-OF} (m, n) := \:$ "$m = \lceil [\psi_1, …, \psi_k]\rceil$ for a valid proof of $\psi_k$ in PA, and $n = \lceil \psi_k \rceil$"
\end{quote}

\noindent Finally, we are ready to define within PA the formula

\begin{quote}
$\phi_{PROVABLE}(n) := \:$ "$\exists m \: \phi_{PROOF-OF}(m,n)$"
\end{quote}

\noindent Importantly, $\phi_{PROVABLE}$ \textbf{does not} represent a decidable function. However, whenever a proof of a formula $\psi$ is provided, the proof can readily be encoded in some $m$, and then $\phi_{PROOF-OF}(m, \lceil \psi \rceil)$ is shown in PA, from which $\phi_{PROVABLE}(\psi)$ follows. This property will suffice to prove the theorem.

\subsubsection{Step 4. Devise a self-referential formula through a variable substitution trick}

The following variable substitution procedure, defined as a function $\mathcal{D}: \mathbb{N} \rightarrow \mathbb{N}$, is a key device in the proof. $\mathcal{D}$ takes as input $\lceil \phi \rceil$, and returns $\lceil \phi^\mathcal{D} \rceil$ where $\phi^{\mathcal{D}} := \phi[x_0: \lceil \phi \rceil]$. The variable $x_0$ is the lexicographically first variable in the alphabet $\Sigma$, and $\phi$ may contain zero or more occurrences of $x_0$. All these occurrences are replaced by $S^{(\lceil \phi \rceil)}0$, and the resulting formula is encoded by $\lceil \rceil$.\footnote{As a reminder, we are abusing notation when we say $\phi[x_0: \lceil \phi \rceil]$: $x_0$ is not replaced by the number $\lceil \phi \rceil$ but instead, by $S...S0$ ($\lceil \phi \rceil$ many $S$'s).}

\begin{quote}
\textbf{Diagonalization:} $\mathcal{D}: \mathbb{N} \rightarrow \mathbb{N}$: Given $n \in \mathbb{N}$,
\begin{itemize}
\item If $\lnot$ $\phi_{FORMULA}(n)$, return $0$.
\item Else, let $n = \lceil \phi \rceil$.
\item Construct the string $s = S^{(n)}0$.
\item Substitute all occurrences of $x_0$ in $\phi$ with $s$, to obtain $\phi^\mathcal{D} := \phi[x_0: s]$. 
\subitem \textit{By abuse of notation, $\phi^{\mathcal{D}} = \phi[x_0: \lceil \phi \rceil]$.}
\item Return $\lceil \phi^{\mathcal{D}} \rceil$.
\end{itemize}
\end{quote}

\noindent The procedure above is recursive, and therefore it is representable in a formula of PA, let's call it $\phi_{DIAG}$, where
$$\vdash \forall y \: (\phi_{DIAG}(n, y) \leftrightarrow y = \mathcal{D}(n))$$
\noindent Now comes a tricky lemma at the heart of the proof.

\begin{quote}
\textbf{Lemma (G\"{o}del's fixed-point lemma or diagonalization lemma):} given a formula $\phi(x)$ in PA, of a single free variable $x$, a sentence $\sigma$ can be constructed such that PA $\vdash \sigma \leftrightarrow \phi(\lceil \sigma \rceil)$.
\end{quote}

\noindent \textbf{Proof:} Consider the formula $\phi_*(x_0) = \exists y \: (\phi_{DIAG}(x_0, y) \land \phi(y))$. This formula is true whenever $x_0 = \lceil \psi \rceil$, for some formula $\psi$ that may have $0$ or more occurrences of variable $x_0$ (those occurrences inside $\psi$ are not to be confused by the outer $x_0 = \lceil \psi \rceil$), $y$ is equal to $\mathcal{D}(x_0) = \mathcal{D}(\lceil \psi \rceil) = \lceil \psi[x_0: \lceil \psi \rceil] \rceil$, and moreover $\phi(y)$ holds where $\phi$ is the starting formula in the Lemma we are seeking to prove. Therefore, $\phi_*(x_0)$ is equivalent to $\phi(\lceil \psi[x_0: \lceil \psi \rceil] \rceil)$.

Now let $n = \lceil \phi_*(x_0) \rceil$. Define $\sigma$, a sentence of no free variables:
$$\sigma \: := \: \exists y \: (\phi_{DIAG}(n, y) \land \phi(y))$$
\noindent Then, by the definition of $\phi_{DIAG}$, we get $\vdash \sigma \leftrightarrow \exists y \: (y = \mathcal{D}(n) \land \phi(y))$.

What is $\mathcal{D}(n)$ in this latter expression? it is $\mathcal{D}(\lceil \phi_*(x_0) \rceil) = \lceil \phi_*[x_0: n] \rceil$. Let's write this down explicitly:
$$\phi_*[x_0: n] = \exists y \: (\phi_{DIAG}(n, y) \land \phi(y)) = \sigma$$
\noindent Therefore, we have get 
$$\vdash \sigma \: \leftrightarrow \: \exists y \: (y = \mathcal{D}(n) \land \phi(y)) \: \leftrightarrow \: \exists y \: (y = \lceil \sigma \rceil \land \phi(y)) \: \leftrightarrow \: \phi(\lceil \sigma \rceil)$$\noindent which is what we wanted to prove.$\square$

In the above proof, it is perhaps tricky at first to see that $\phi_{DIAG}(\lceil \exists y \: (\phi_{DIAG}(x_0, y) \land \phi(y)) \rceil, y)$ assigns to $y$ precisely $\lceil \sigma \rceil$. Figure 1 makes this clear.

\begin{figure}
\centering
\includegraphics[width=320pt]{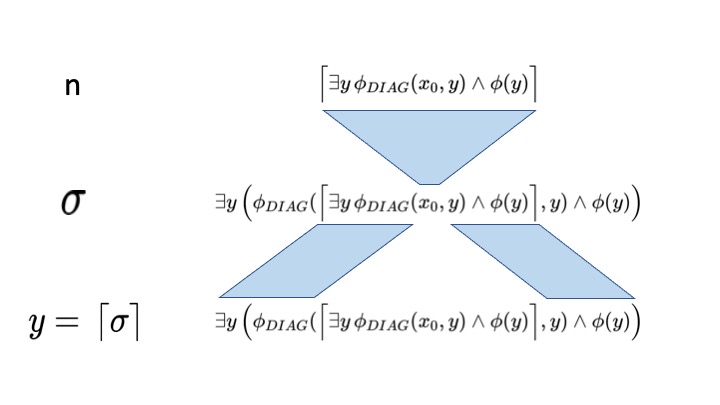}
\caption{\textbf{Construction of sentence $\sigma$ in the diagonalization lemma}. The formula $\phi_{DIAG}$ implements diagonalization: given $n$ and $y$, $\phi_{DIAG}$ asserts that $y$ is the encoding of a formula obtained as follows. First, assert that $n = \lceil \phi \rceil$. Next, find all occurrences of $x_0$ within $\phi$ and replace them with $S^{(n)}$. Finally, encode the result in an integer and assert $y$ is that integer. Now, if we construct $\sigma$ as depicted, $\phi_{DIAG}$ will replace $x_0$ with $n$ and the result will be - voila! - $\sigma$ again so $y$ becomes $\lceil \sigma \rceil$.}
\label{fig: diag}
\end{figure}

\subsubsection{Step 5: conclude the proof}
Therefore, given any formula $\phi$ of a single free variable, we can construct a sentence $\sigma$ of no free variables, such that PA $\vdash \sigma \leftrightarrow \phi(\lceil \sigma \rceil)$. This is a powerful lemma! The rest of the proof is easy.

First, an important result follows immediately: Tarski’s inexistence of a truth definition.

\begin{quote}
\textbf{Tarski's inexistence of Truth definition.} Imagine that we had a table, or a computable function $Truth: \mathbb{N} \rightarrow \mathbb{N}$, such that given any sentence $\tau$ as input, $Truth(\lceil \tau \rceil)$ returns the truth value of $\tau$.

Use the above lemma with $\phi(x) = \lnot Truth(x)$ to construct $\sigma$ such that $\sigma \leftrightarrow \phi(\lceil \sigma \rceil) \leftrightarrow \lnot Truth(\lceil \sigma \rceil)$. 

Is $\sigma$ true? Well, if it is, then $\lnot Truth(\lceil \sigma \rceil)$ is true, in which case $Truth(\lceil \sigma \rceil)$ is false, therefore $\sigma$ is false. And vice versa. We conclude that $Truth(x)$ cannot be a computable function.$\square$
\end{quote}

\noindent Now let’s turn to the incompleteness theorem. Assume PA is complete so that given any true sentence $\tau$, there is a PA proof of $\tau$. We use the lemma above with $\phi = \lnot \phi_{PROVABLE}(x)$. Let $\sigma$ be such that $\vdash \sigma \leftrightarrow \lnot \phi_{PROVABLE}(\lceil \sigma \rceil)$. Either $\sigma$ or $\lnot \sigma$ is true and therefore provable in PA.

\underline{Case1.} $\sigma$ is provable. Let $[\psi_1, ..., \psi_k = \sigma]$ be a proof, and let $m^* = \lceil [\psi_1, ..., \psi_k] \rceil$. By $\vdash \sigma \leftrightarrow \lnot \phi_{PROVABLE}(\lceil \sigma \rceil)$ we get $\vdash \lnot \exists m \: \phi_{PROOF-OF}(m, \lceil \sigma \rceil)$ and simultaneously we have $\vdash \phi_{PROOF-OF}(m^*, \lceil \sigma \rceil)$, therefore PA is inconsistent.

\underline{Case 2.} $\lnot \sigma$ is provable. Then $\phi_{PROVABLE}(\lceil \sigma \rceil)$ $=$ $\exists m \:\phi_{PROOF-OF}(m, \lceil \sigma \rceil)$ follows. From the proof of $\lnot \sigma$, it follows that for any given $m \in \mathbb{N}$, $\lnot \phi_{PROOF-OF}(m, \lceil \sigma \rceil)$ otherwise PA would derive $\sigma$, a contradiction. Here, in his original proof G\"{o}del assumed $\omega$-consistency of arithmetic: 
\begin{quote}
\textbf{$\omega$-consistency: } there is no formula $\psi(x)$ such that for every $m$, $\vdash \lnot \psi(m)$ and also $\vdash \exists m \: \psi(m)$. 
\end{quote}

\noindent We now have PA $\vdash \lnot \sigma$ therefore for every $m$, $\vdash \lnot \phi_{PROOF-OF}(m, \lceil \sigma \rceil)$, and also $\exists m \:\phi_{PROOF-OF}(m, \lceil \sigma \rceil)$, contradicting $\omega$-consistency and concluding the proof.$\square$

However, $\omega$-consistency is a strong assumption. For example, we could augment PA with the axiom $\lnot \sigma$. The resulting theory would be $\omega$-inconsistent but still not inconsistent in the sense of proving a contradiction.

\subsubsection{Rosser's trick and the G\"{o}del-Rosser Incompleteness Theorem}
Fortunately, we can fix the above proof with what is known as Rosser's trick, a simple modification to the formula $\phi_{PROVABLE}$ that makes the theorem a lot more powerful.

We alter the definition of $\phi_{PROOF-OF}(m, n)$. First, let $\phi_{NOT}(m, n) := \:$ $\exists \psi \: m = \lceil \psi \rceil \land n = \lceil \lnot \psi \rceil$.

\begin{quote}
\textbf{Rosser's trick:}
$$\phi^R_{PROOF-OF}(m, n) := \: \phi_{PROOF-OF}(m, n) \land \lnot (\exists k \leq m \: (\exists l \:(\phi_{NOT}(n, l) \land \phi_{PROOF-OF}(k, l)))$$
\end{quote}

\noindent Basically, $\phi^R_{PROOF-OF}(m, n)$ says that $m$ is the encoding of a proof of a formula encoded by $n$, and there is no $k \leq m$ that is the encoding of a proof of the negation of the formula encoded by $n$. Intuitively,

\begin{quote}
$\phi^R_{PROOF-OF}(m, n) :=$ "$m$ encodes a proof of the formula encoded by $n$, and there is no shorter proof of this formula's negation"
\end{quote}

\noindent Now we can proceed as before and define $\phi^R_{PROVABLE}(n) = \exists m \:\phi^R_{PROOF-OF}(m, n)$, and use the lemma with $\phi^R = \lnot \phi^R_{PROVABLE}(x)$. Let $\sigma^R$ be such that $\vdash \sigma^R \: \leftrightarrow \: \lnot \phi^R_{PROVABLE}(\lceil \sigma^R \rceil)$.

Let's go back to the problematic case where $\lnot \sigma$ was provable, which required $\omega$-consistency to get a contradiction. Let $\lnot \sigma^R$ be provable and let $[\phi_1, ..., \phi_k = \lnot \sigma^R]$ be a proof. If PA is consistent, this implies there is no shorter proof of $\sigma^R$, something that can be checked in a primitive recursive manner, therefore $\vdash \phi_{PROVABLE}^R(\lceil \lnot \sigma^R \rceil)$. By diagonalization we also have $\vdash \phi_{PROVABLE}^R(\lceil \sigma^R \rceil)$, which asserts that there is no shorter proof of $\lnot \sigma^R$. Both statements cannot be true simultaneously: $\vdash \phi_{PROVABLE}^R(\lceil \sigma^R \rceil)$ can be checked to be false by searching for any proof of $\sigma^R$ that is shorter than the proof of $\lnot \sigma^R$.$\square$

\subsubsection{A two-sentence liar's paradox}

The incompleteness theorem is a version of the liar's paradox "this sentence is not true", modified to "this statement is not provable". Other versions of the liar's paradox exist, such as the two-sentence version:
\begin{quote}
A: The following statement is false.

B: The preceding statement is true.
\end{quote}

\noindent Can we devise an incompleteness proof based on this version? The following lemma constructs two sentences that "talk" about each other.\footnote{When I  wrote this section and the next, I thought the results were new and was not aware of previous generalizations of diagonalization (see Buldt 2014). Still, my constructions are clean ways to generalize diagonalization so I keep them here.}

\noindent \textbf{Lemma: Two-Sentence Liar's Paradox.} Given a formula $\phi(x)$ in PA, of a single free variable $x$, sentences $\sigma$ and $\tau$ can be constructed such that  PA$\vdash \sigma \leftrightarrow \phi(\lceil \tau \rceil)$ and PA $\vdash \tau \leftrightarrow \lnot \phi(\lceil \sigma \rceil)$.

\noindent \textbf{Proof:} First, we define a function $\mathcal{R}: \mathbb{N} \rightarrow \mathbb{N}$ that will play a role analogous to $\mathcal{D}$ of the diagonalization lemma.

\begin{quote}
\textbf{Reversal:} $\mathcal{R}: \mathbb{N} \rightarrow \mathbb{N}$: Given $n \in \mathbb{N}$,
\begin{enumerate}
\item If $\lnot$ $\phi_{FORMULA}(n)$, return $0$; else, let $n = \lceil \phi \rceil$.
\item If $\phi \not = \exists y \: \psi(x_0, y) \land \chi(y)$, where $\psi$ has at most $x_0, y$ free, $x_0$ is the specific first variable in $\Sigma$, and $\chi$ has at most $y$ free, return $0$.
\item If $\chi = \lnot \chi'$, then let $\chi^{\mathcal{R}} = \chi'$ else let $\chi^{\mathcal{R}} = \lnot \chi$.
\item Let $n' = \lceil \exists y \: \psi(x_0, y) \land \chi^{\mathcal{R}} \rceil$.
\item Construct the string $s = S^{(n')}0$.
\item Substitute all occurrences of $x_0$ in $\psi$ with $s$, to obtain $\psi^\mathcal{R} := \psi[x_0: s]$. 
\item Construct the formula $\phi^{\mathcal{R}} = \exists y \: \psi^{\mathcal{R}}(y)\land \chi^{\mathcal{R}}(y)$
\item Return $\lceil \phi^{\mathcal{R}} \rceil$.
\end{enumerate}
\end{quote}

\noindent The procedure above is (primitive) recursive, and therefore is representable with a formula in PA, let's call it $\phi_{REV}$, where
$$\vdash \forall y \: (\phi_{REV}(n, y) \leftrightarrow y = \mathcal{R}(n))$$
\noindent Starting with $\phi$ in the lemma statement, let $y$ be a variable of $\Sigma$ not appearing within $\phi$ or $\phi_{REV}$. Let $n = \lceil \exists y \: \phi_{REV}(x_0, y) \land \phi(y) \rceil$ and let $m = \lceil \exists y \: \phi_{REV}(x_0, y) \land \lnot \phi(y) \rceil$, unless $\phi = \lnot \phi'$ in which case we let $m = \lceil \exists y \: \phi_{REV}(x_0, y) \land \phi'(y) \rceil$. Let $\sigma = \exists y \: \phi_{REV}(n, y) \land \phi(y)$ and $\tau = \exists y \: \phi_{REV}(m, y) \land \lnot \phi(y)$, unless $\phi = \lnot \phi'$, in which case we let $\tau = \exists y \: \phi_{REV}(m, y) \land \phi'(y)$. 

Calculating $\sigma$, we get:
$$\sigma = \exists y \: \phi_{REV}(\lceil \exists y \: \phi_{REV}(x_0, y) \land \phi(y) \rceil, y) \land \phi(y)  \leftrightarrow$$\noindent (the outer $\phi_{REV}(n, y)$ will set $n' = \lceil \exists y \: \phi_{REV}(x_0, y) \land \lnot \phi(y) \rceil$ in step 5, and set $y$ to $\lceil \phi^{\mathcal{R}} \rceil = \lceil \exists y \: \phi_{REV}(n', y) \land \lnot \phi(y) \rceil$).
$$\exists y \: (y = \lceil \exists y \: \phi_{REV}(\lceil \exists y \: \phi_{REV}(x_0, y) \land \lnot \phi(y) \rceil, y) \land \lnot \phi(y) \rceil \land \phi(y) \leftrightarrow$$
$$\phi(\lceil \exists y \: \phi_{REV}(\lceil \exists y \: \phi_{REV}(x_0, y) \land \lnot \phi(y) \rceil), y) \land \lnot \phi(y) \rceil) \leftrightarrow \phi(\lceil \tau \rceil)
$$
\noindent Then, calculating $\tau$ we get:
$$\tau = \exists y \: \phi_{REV}(m, y) \land \lnot \phi(y) = \exists y \: \phi_{REV}(\lceil \exists y \: \phi_{REV}(x_0, y) \land \lnot \phi(y) \rceil, y) \land \lnot \phi(y) \leftrightarrow$$
$$\exists y\: (y = \lceil \exists y \: \phi_{REV}(\lceil \exists y \: \phi_{REV}(x_0, y) \land \phi(y) \rceil, y) \land \phi(y) \rceil) \land \lnot \phi(y) \leftrightarrow$$
$$\lnot \phi(\lceil \exists y \: \phi_{REV}(\lceil \exists y \: \phi_{REV}(x_0, y) \land \phi(y) \rceil, y) \land \phi(y) \rceil) \leftrightarrow \lnot \phi(\lceil \sigma \rceil)
$$

\noindent This concludes the proof of the lemma.$\square$

Now the  theorem follows from the two-sentence lemma by letting $\phi(x) = \phi_{PROVABLE}^{\mathcal{R}} (x)$. We construct $\sigma, \tau$ such that $\vdash \sigma \leftrightarrow \phi_{PROVABLE}^{\mathcal{R}}(\lceil \tau \rceil)$ and $\vdash \tau \leftrightarrow \lnot \phi_{PROVABLE}^{\mathcal{R}}(\lceil \sigma \rceil)$. Let $\vdash \sigma$. Then $\vdash \phi_{PROVABLE}^{\mathcal{R}}(\lceil \sigma \rceil)$, therefore  $\vdash \lnot \tau$ and simultaneously $\phi_{PROVABLE}^{\mathcal{R}}(\lceil \tau \rceil)$, inducing a contradiction. Otherwise, let $\vdash \lnot \sigma$: then $\vdash \lnot \phi_{PROVABLE}^{\mathcal{R}}(\lceil \tau \rceil)$. If $\vdash \tau$ we get a contradiction, whereas if $\vdash \lnot \tau$ we get $\phi_{PROVABLE}^{\mathcal{R}}(\lceil \sigma \rceil)$, a contradiction. Therefore, neither $\sigma$ nor $\lnot \sigma$ is provable.$\square$

\subsubsection{Further generalization to multi-sentence liar's paradox }

There are multi-sentence versions of the liar's paradox, such as circular constructions of the form:
\begin{quote}
$A_1$: Sentence $A_2$ is false.

$A_2$: Sentence $A_3$ is false.

...

$A_k$: Sentence $A_1$ is false.
\end{quote}

\noindent This leads to a paradox whenever $k$ is odd: if $A_1$ is true then $A_k$ is true, therefore $A_1$ is false and so on. We can further generalize the diagonalization lemma to cover these cases and more. We will construct sentences $\sigma_1, ..., \sigma_k$, whereby each $\sigma_i \leftrightarrow \psi_i(\lceil \sigma_{f(i)} \rceil)$ is a theorem of PA, where $f$ is any function on $\{1, ..., k\}$ and the $\psi_i$s are any formulas.

\noindent \textbf{Lemma: Generalized Diagonalization.} Given single-variable formulas $\psi_1, ..., \psi_k$ in PA and function $f: \{1, ..., k\} \rightarrow \{1, ..., k\}$, sentences $\sigma_1, ..., \sigma_k$ can be constructed such that $\forall_{1 \leq i \leq k} \:$ PA $\vdash \sigma_i \leftrightarrow \psi_i(\lceil \sigma_{f(i)} \rceil)$.

\noindent \textbf{Proof:} Fix single-variable formulas $\psi_1, ..., \psi_k$ and function $f: \{1, ..., k\} \rightarrow \{1, ..., k\}$. For example we could have $f(i) := i+1$ mod $k$ and $\psi_i(y) = \lnot \phi(y)$, to produce a circular liar's paradox whenever $k$ is odd.

First we define $\mathcal{GD}: \mathbb{N} \rightarrow \mathbb{N}$, a generalized version of the diagonalization function.

\begin{quote}
\textbf{Generalized Diagonalization:} $\mathcal{GD}: \mathbb{N} \rightarrow \mathbb{N}$: Given $n \in \mathbb{N}$,
\begin{enumerate}
\item If $n \not = \lceil [\exists y \: \phi_{GD}(x_0, y) \land \psi_i(y), i, \phi(y)] \rceil$, return $0$.
\item Let $n' = \lceil [\exists y \: \phi_{GD}(x_0, y) \land \psi_{f(i)}, f(i), \phi(y)] \rceil$.
\item Construct the string $s = S^{(n')}0$.
\item Substitute all occurrences of $x_0$ in $\phi_{GD}$ with $s$, to obtain $\phi_{GD}(s, y)$.
\item Return $\lceil \phi_{GD}(s, y) \land \psi_{f(i)}  \rceil$.
\end{enumerate}
\end{quote}

\noindent The procedure above is (primitive) recursive, and therefore is representable with a formula in PA, let's call it $\phi_{GD}$, where
$$\vdash \forall y \: (\phi_{GD}(n, y) \leftrightarrow y = \mathcal{GD}(n))$$For $1 \leq i \leq k$ let $n_i = \lceil [\exists y \: \phi_{GD}(x_0, y) \land \psi_i(y), i ] \rceil$, and let $\sigma_i = \exists y \: \phi_{GD}(n_i, y) \land \psi_i(y)$, where $y$ is a variable not appearing within $\phi_{GD}$ or any $\psi_i$. Then, by calculating the return value of $\mathcal{GD}(n_i)$, we get:
$$\sigma_i \leftrightarrow \exists y \: \phi_{GD}(\lceil [\exists y \: \phi_{GD}(x_0, y) \land \psi_i(y), i ] \rceil, y) \land \psi_i(y) \leftrightarrow$$ $$(y = \lceil \exists y \: \phi_{GD}(\lceil [\exists y \: \phi_{GD}(x_0, y) \land \psi_{f(i)}, f(i)] \rceil, y) \land \psi_{f(i)}(y) \rceil) \land \psi_i(y) \leftrightarrow$$ $$(y = \lceil \sigma_{f(i)} \rceil \land \psi_i(y)) \: \leftrightarrow \: \psi_i(\lceil \sigma_{f(i)} \rceil)$$
\noindent This concludes the proof of the generalized lemma. $\square$

\subsubsection{Direct self-reference}

Let's examine the statement $\sigma$ constructed by the diagonalization lemma: 
$$\sigma \: := \: \exists y \: \phi_{DIAG}(\lceil \exists y \: \phi_{DIAG}(x_0, y) \land \lnot \phi_{PROVABLE}(y) \rceil, y) \land \lnot \phi_{PROVABLE}(y)
$$
In English, "there is a number $y$ such that there is a formula encoded by $y$, which is the result of substituting within the formula encoded by the number $\lceil$[\textit{there is a number $y$ such that there is a formula encoded by $y$, which is the result of substituting within the formula encoded by the number $x_0$ all occurrences of the character `$x_0$' with the number $x_0$,}\footnote{The occurrences of $x_0$ within the "formula encoded by (outer) $x_0$" should not be confused with the outer occurrence, which is bound. This is an important point, so an example is in order. Say $x_0 = \lceil x_0 = SS0 \rceil$. This is perfectly legit. The outer $x_0$ is bound to a number, specifically the encoding of the formula $x_0 = SS0$, whereby the inner $x_0$ is simply a letter of $\Sigma$ within that formula. Then, when we are substituting within the formula encoded by $x_0$ all occurrences of `$x_0$' with the number $x_0$, we are substituting within "$x_0 = SS0$" all occurrences of `$x_0$' with the number $\lceil x_0 = SS0 \rceil$, so the inner formula turns into $S^{(\lceil x_0 = SS0 \rceil)}0 = SS0$. Kripke uses special terminology for this situation: the outer occurrence of $x_0$ is "used" and the inner occurrences are "mentioned".} \textit{and this resulting formula is not provable}]$\rceil$ all occurrences of the character `$x_0$' with the number $\lceil$[\textit{there is a number $y$ such that there is a formula encoded by $y$, which is the result of substituting within the formula encoded by the number $x_0$ all occurrences of the character `$x_0$' with the number $x_0$, and this resulting formula is not provable}]$\rceil$, and this resulting formula is not provable".

Whether the above paragraph is English is debatable. What is clear is that the formula is not directly self-referential: it employs a laborious trick of variable substitution (represented by $\phi_{DIAG}$) to be able to refer to itself. Is there a way for a formula to directly refer to itself without some faulty circularity?

In a short, elegant paper, Saul Kripke describes a few different ways for deploying direct reference (Kripke 2021). Kripke attributes the following way to Raymond Smullyan through personal communication (no reference provided):

First, let $x_0$ be the lexicographically first variable of $\Sigma$. Let $A_1(x_0), A_2(x_0), ...$ be an enumeration of all formulas that contain no free variables other than $x_0$. Let the original encoding $\lceil \rceil$ be such that the smallest prime number employed is $2$, resulting in G\"{o}del numbers that are always odd. Define $\lceil \rceil^*$ to be a new encoding that coincides with $\lceil \rceil$ except as follows. For each $n$, let $k_n = \lceil \exists x_0 \: (x_0 = n \land A_n(x_0)) \rceil$. Now the following formula gets a special encoding:
$$\lceil \exists x_0 \: (x_0 = 2k_n \land A_n(x_0)) \rceil^* := 2k_n$$
\noindent In this manner, every formula $A_n(x_0)$ has an instance as above, which gets the encoding $2k_n$, asserting that its own encoding satisfies $A_n(2k_n)$. Then, the theorem proceeds as before with $A_n(x_0)$ being $\lnot \phi_{PROVABLE}(x_0)$ or the Rosser version $\lnot \phi_{PROVABLE}^R(x_0)$.

\section{Other ways to prove the incompleteness theorem}

I outline informally two other strategies to prove G\"{o}del's Incompleteness theorem: through Kolmogorov complexity and through Turing machines.

\noindent\textbf{Kolmogorov complexity.} The Kolmogorov complexity $K(n)$ of a number $n \in \mathbb{N}$ is the shortest computer program that outputs $n$ and terminates. The language can be fixed - say C or Python - and the program length can be in characters or in bits. For every $L \in \mathbb{N}$, there are numbers $n \in \mathbb{N}$ such that $K(n) > L$. The simple way to see this is that there are at most $2^L$ programs of length $L$. Therefore, some numbers $> 2^L$ are not the output of a program of length $\leq L$. Let's fix an axiomatic language $\mathcal{L}$ that is rich enough to represent proofs of statements "$K(n) > m$". Consider now the following program $PRINT(L)$ := "go through all proofs of $\mathcal{L}$ lexicographically, find the first number $n$ proven to have $K(n) > L$, print $n$ and terminate". Let's assume that the length of $PRINT()$ is $L_{PRINT}$ bits. Now let's run $PRINT(L_{PRINT})$. It will run until it finds $n$ with $K(n) > L_{PRINT}$, print $n$ and terminate. Therefore $K(n) \leq L_{PRINT}$, a contradiction. The conclusion is that $PRINT(L_{PRINT})$ will loop forever, and consequently $\mathcal{L}$ has no proof of the form $K(n) > L_{PRINT}$. This is Chaitin's theorem: for any rich enough system $\mathcal{L}$ there is an $L$ such that no $n$ can be proven to have $K(n) > L$.

\noindent\textbf{Turing Machines.} The proof of incompleteness through the halting problem in Turing machines is well known. A proof that incorporates the Rosser improvement of not requiring $\omega$-consistency is sketched by Scott Aaronson in his superb blog [1]. First, define the following variant of the halting problem:

\begin{quote}
\textbf{The Consistent Guessing Problem.} Given a description of a Turing machine $M$,
\begin{enumerate}
\item If $M$ accepts on a blank tape, accept.
\item If $M$ rejects on a blank tape, reject.
\item If $M$ loops forever, accept and halt.
\end{enumerate}
\end{quote}

\noindent It is easy to argue that there is no Turing machine that solves this problem: let $P$ be such a machine. Modify $P$ to obtain $Q$: "given $M$, if $M(M)$\footnote{The standard notation is $M(\langle M \rangle)$ to denote the description of a machine $M$ by $\langle M \rangle$. There is no confusion in abusing notation with $M(M)$: whenever a machine is an input to a machine, it turns into its syntactically valid description.} accepts, reject; if $M(M)$ rejects, accept; otherwise accept and halt". Now ask what does $Q(Q)$ yield? Well, reading from the description above, "given  $M=Q$, if $Q(Q)$ accepts, reject; if $Q(Q)$ rejects, accept..." Now we are in trouble: on input $Q$, $Q$ will do the opposite of what $Q$ does on input $Q$, contradiction. Therefore $P$ does not exist.

Now assume $\mathcal{L}$ is rich enough to formalize Turing machines and is complete and consistent. Then, solve the Consistent Guessing Problem above as follows: given $M$, enumerate lexicographically all possible proofs of $\mathcal{L}$ until a proof of "$M$ rejects on a blank tape" or "$M$ does not reject on a blank tape" appears; if the former, reject and if the latter, accept. Therefore, such a complete and consistent $\mathcal{L}$ is impossible.

Both of these strategies are excellent for introducing the theorem to computational people. However, they both mask the complicated task of formalizing computation in axiomatic systems.

\section{G\"{o}del's Second Incompleteness Theorem}

G\"{o}del's Second Incompleteness theorem informally states that a rich enough axiomatic system cannot prove its own consistency. This theorem follows easily from a combination of the first theorem and the following property of $\phi_{PROVABLE}$ (which in this section we abbreviate to $\phi_{PR}$): $\vdash \phi_{PR}(\lceil \phi \rceil) \rightarrow \phi_{PR}(\phi_{PR}(\lceil \phi \rceil))$. The latter property is technically hard to prove and beyond the scope of this manuscript. For a complete exposition see for instance Swierczkowski (2003).\footnote{It is technically messy to show that for any $\Sigma$-formula $\phi$, if all free variables in its encoding $\lceil \phi \rceil$ are handled properly, then $\vdash \phi \rightarrow \phi_{PR}(\lceil \phi \rceil)$. Because $\phi_{PR}$ is a $\Sigma$ formula, the result then follows. As an example of disasters that  occur when variables are not properly handled, consider the $\Sigma$-formula $x_0 = x_1$. If it were the case that $\vdash x_0 = x_1 \rightarrow \phi_{PR}(\lceil x_0 = x_1 \rceil)$ then since $x_0, x_1$ do not occur in $\phi_{PR}(\lceil x_0 = x_1 \rceil)$ (the latter being simply $\phi_{PR}(n)$ for a specific $n \in \mathbb{N}$ with no free variables), it would follow that $\vdash \phi_{PR}(\lceil x_0 = x_1 \rceil)$ (simply replace $x_0$ with $0$ and $x_1$ with $0$ to obtain $true \rightarrow \phi_{PR}(\lceil x_0 = x_1 \rceil)$). Finally, from $\vdash \phi_{PR}(\lceil x_0 = x_1 \rceil)$ we conclude $x_0 = x_1$, a contradiction if we now let $x_0 = S0, x_1 = 0$.}

More formally, if $\mathcal{L}$ is rich enough to express provability, $\phi_{PR}$, then $\mathcal{L} \not \vdash \lnot \phi_{PR}(S0 = 0)$. Consistency is expressed with the notion that the absurdity $1 = 0$ is not provable in the system. The proof is by contradiction: assuming that $\mathcal{L}$ proves $\lnot \phi_{PR}(S0 = 0)$, it follows that $\mathcal{L}$ also proves $S0 = 0$. Therefore either $\mathcal{L}$ is inconsistent, or it cannot prove $\lnot \phi_{PR}(S0 = 0)$.

The truth predicate $\phi_{PR}$ satisfies the following conditions, known as the Hilbert-Bernays conditions, which are what is needed for the second incompleteness proof to go through:

\begin{quote}
\textbf{The Hilbert-Bernays provability conditions:}
\begin{enumerate}[label=(\roman*)]
\item If $\mathcal{L} \vdash \phi$ then $\mathcal{L}  \vdash \phi_{PR}(\lceil \phi \rceil)$.
\item $\mathcal{L}  \vdash \phi_{PR}(\lceil \phi \rceil) \rightarrow \phi_{PR}(\lceil \phi_{PR}(\lceil \phi \rceil) \rceil)$.
\item $\mathcal{L}  \vdash \phi_{PR}(\lceil \phi \rceil) \land \phi_{PR}(\lceil \phi \rightarrow \psi \rceil) \rightarrow \phi_{PR}(\lceil \psi \rceil)$.
\end{enumerate}
\end{quote}

\noindent\textbf{Theorem (G\"{o}del's Second Incompleteness Theorem.} $\mathcal{L}  \not \vdash \lnot \phi_{PR}(S0 = 0)$.

\noindent In the proof below, all derivations are within $\mathcal{L}$ and $\mathcal{L} \vdash$ in front of every derived formula is implied.

\noindent \textbf{Proof:} The goal is to show that $\lnot \phi_{PR}(S0=0)$ leads to a contradiction. We might as well start where the first theorem left off: using  diagonalization, let $\sigma$ be a sentence such that $\sigma \leftrightarrow \lnot \phi_{PR}(\lceil \sigma \rceil)$. Deriving $\sigma$ (or $\lnot \sigma)$ leads to a contradiction, so it suffices to show  $\lnot \phi_{PR}(S0=0) \rightarrow \sigma$. 

This will be accomplished by deriving $\phi_{PR}(\lceil \sigma \rceil) \rightarrow \phi_{PR}(\lceil \lnot \sigma \rceil)$ as an intermediate step: starting from the diagonalization lemma, $\phi_{PR}(\lceil \sigma \rceil) \rightarrow \lnot \sigma$, by (i) $\phi_{PR}(\phi_{PR}(\lceil \sigma \rceil) \rightarrow \lnot \sigma)$. By (iii), this yields $\phi_{PR}(\phi_{PR}(\lceil \sigma \rceil)) \rightarrow \phi_{PR}(\lceil \lnot \sigma \rceil)$. Using (ii) this yields $\phi_{PR}(\lceil \sigma \rceil) \rightarrow \phi_{PR}(\lceil \lnot \sigma \rceil)$ (*).

Then $\phi_{PR}(\lceil \sigma \rceil) \rightarrow (\phi_{PR}(\lceil \lnot \sigma \rceil) \land \phi_{PR}(\lceil \sigma \rceil)$ follows, and by definition of $\phi_{PR}$ we have $\phi_{PR}(\lceil \phi \rceil) \land \phi_{PR}(\lceil \psi \rceil) \rightarrow \phi_{PR}(\lceil \phi \land \psi \rceil)$ from which we get $\phi_{PR}(\lceil \sigma \rceil) \rightarrow (\phi_{PR}(\lceil \lnot \sigma \land \sigma \rceil)$ which is equivalent to $\phi_{PR}(\lceil \sigma \rceil) \rightarrow \phi_{PR}(S0 = 0)$.

Reversing the direction of implication, $\lnot \phi_{PR}(S0 = 0) \rightarrow \lnot \phi_{PR}(\lceil \sigma \rceil)$, which by $\sigma \leftrightarrow \lnot \phi_{PR}(\sigma)$ yields $\lnot \phi_{PR}(S0=0) \rightarrow \sigma$.$\square$

\section{Incompleteness and Completeness}

Gödel's incompleteness theorems introduce the notion of a proposition that is \textit{independent} of a language: the statement $\sigma$ cannot be proved or disproved, and is called \textit{independent} of PA. Another famous theorem courtesy of Gödel is the Completeness theorem for axiomatic systems of first-order logic with equality - roughly, systems like PA with domain-specific axioms, logical axioms, equality axioms, and quantifiers $\forall, \exists$ whose variables range over "first-order" objects: the members of the domain.\footnote{Quantifiers cannot range over second-order or higher objects, such as properties of members of the domain. Discussing second-order logic is beyond the scope here.} Informally, the theorem states that any \textit{necessarily} true statement in a first-order theory is provable.

At first sight, completeness and incompleteness seem contradictory. Discussing completeness of first-order logic is beyond the scope of this manuscript. Instead, I will briefly outline the notions necessary to address and clarify the seeming contradiction. 

\textbf{Necessarily true} involves the notion of a \textit{model} for the axiomatic theory. Briefly, a model $\mathcal{M}$ consists of a \textbf{domain} for the terms and an \textbf{interpretation} of the function and relation symbols of the theory, \textbf{such that all axioms are true} under that interpretation. In the case of PA the standard model is $\mathbb{N}$ with standard interpretation of $S, <, >, =, +, *$. However, this is not the only possible model.\footnote{In fact, the "upward" Löwenheim–Skolem theorem states that for every enumerable language $\mathcal{L}$, for every model $\mathcal{M}$ of $\mathcal{L}$ of cardinality $\kappa$, and for every cardinality $\lambda > \kappa$, we can obtain a model $\mathcal{N}\supset \mathcal{M}$ of cardinality $\lambda$ such that the truth of any proposition $\phi$ according to the two models agrees, as long as all variables of $\phi$ range over the smaller model. This means that starting from $\mathbb{N}$ we can extend to models of PA of arbitrarily large cardinality. A theory is \textbf{categorical} if it has only one model up to isomorphism. The Löwenheim–Skolem theorem implies that a first-order theory with an infinite model, such as PA with $\mathbb{N}$, can never be categorical.} For instance, the domain could contain additional elements that are not in the successor chain of $0$. The Completeness theorem says that given a consistent first-order theory $\mathcal{L}$, a sentence $\tau$ is provable iff $\tau$ is true in every possible model of $\mathcal{L}$. Consequently, if a sentence such as $\sigma$ is neither provable nor refutable in PA, then there exist models $\mathcal{M}$ and $\mathcal{M}'$ with $\sigma$  true in $\mathcal{M}$ and $\lnot \sigma$  true in $\mathcal{M}'$. Moreover, we can extend PA to PA$+\{ \sigma\}$ and to PA$+\{ \lnot \sigma\}$ corresponding to these two models, respectively, and both resulting theories are consistent assuming that PA is consistent.

How can this be possible? Let's take a look at the two extensions. In PA$+\{ \sigma\}$, we would have $\vdash \sigma$ and $\vdash \sigma \leftrightarrow \lnot \phi_{PROVABLE}(\sigma)$ therefore we get $\vdash \lnot \phi_{PROVABLE}(\sigma)$ which seems to contradict $\vdash \sigma$. It actually doesn't. Recall that $\lnot \phi_{PROVABLE}(\sigma) \: := \: \lnot \exists m \: \phi_{PROOF-OF}(m, \sigma)$ states that $\sigma$ is not provable using PA's axioms, and does not state that $\sigma$ is not provable in PA$+\{ \sigma \}$. Hence there is no contradiction.

OK, how about PA$+\{ \lnot \sigma\}$? Then we get $\vdash \lnot \sigma$ and $\vdash \phi_{PROVABLE}(\sigma)$, i.e., $\vdash \exists m \: \phi_{PROOF-OF}(m, \sigma)$. Wouldn't that be a contradiction? Not really, because the resulting theory would be $\omega$-inconsistent but still consistent. How about the Rosser version of the theorem, leading to $\vdash \exists m \: \phi_{PROOF-OF}^R(m, \sigma)$? Recall that now $m$ is asserted to be smaller than any $n$ such as $\phi_{PROOF-OF}^R(n, \lnot \sigma)$, and we do have $n = \lceil [\lnot \sigma ]\rceil$ because $\lnot \sigma$ is now an axiom therefore $[\lnot \sigma ]$ is a proof. Here is where the nonstandard nature of the model $\mathcal{M}'$ shows up: $m$ will be a number not in the chain of successors of $0$. Such an $m$ would satisfy $\phi_{PROOF-OF}^R(m, \sigma)$ without encoding an actual proof of $\sigma$. Generally, $m$ will be a new constant in the domain postulated to satisfy all the required properties. For more details, reading the Henkin proof of the Completeness theorem would be instructive (Arjona and Alonso 2014).

Actually, any consistent first-order theory $\mathcal{L}$ can be extended to a complete theory $\mathcal{L}'$ where every sentence $\tau$ is either provable or disprovable! This is simple to show: (1) lexicographically order all sentences, $\sigma_1, \sigma_2, ...$ of $\mathcal{L}$; (2) find the first independent sentence $\not \vdash \sigma_i$ and $\not \vdash \lnot \sigma_i$, and define this to be $\tau_1 := \sigma_i$; (3) let $\mathcal{L}_1$ be $\mathcal{L} \cup \tau_1$,  extending $\mathcal{L}$ with $\tau_1$ as axiom; (4) proceed to define $\tau_2, \mathcal{L}_2$, $\tau_3, \mathcal{L}_3, ... \tau_n, \mathcal{L}_n, ...$ for all $n$; (5) let $\mathcal{L}' = \cup{\mathcal{L}_i}$. It is easy to show that $\mathcal{L}'$ is consistent and complete, assuming $\mathcal{L}$ is consistent. How come this does not contradict the incompleteness theorem? The key observation here is that $\{\tau_1, \tau_2, ...\}$ is not a recursive collection of sentences. Because of that, $\phi_{PROVABLE}$ cannot be written down: provability cannot be represented in $\mathcal{L}'$. 

The upshot is that in the standard model $\mathbb{N}$, $\sigma$ is true but not provable in PA. However, due to the Completeness theorem we can extend PA either with $\sigma$ or with $\lnot \sigma$ and maintain consistency. In the first case the new theory is consistent because the $\phi_{PROVABLE}$ predicate does not include the new axiom $\sigma$. In the second case we get $\omega$-inconsistency but consistency, or in the Rosser version we necessarily get an extended domain beyond $\mathbb{N}$ that includes new non-numerical constants that make the statement $\phi_{PROVABLE}(\sigma)$ true without a numerical encoding of an actual proof being produced. Finally, theoretically PA can be extended to a complete theory (in infinitely many ways, actually), but the new infinite collection of axioms will not be computable. In this case the incompleteness theorem does not apply, however such a theory is useless because it can never be written down.

\section{Discussion}
The broader implications of Gödel's Incompleteness theorems on foundational mathematics, logic, philosophy, computer science, physics, and the nature of the human mind, have been debated extensively. Here, I state some of my personal positions as succinctly as I can.

\textbf{Mathematics.} Any mathematical argument can be expressed as a finite number of presuppositions that lead to a conclusion through a number of first-order logical inferences that can be machine-verified syntactically (see also Kripke 2013).\footnote{Every mathematical proof is finite, therefore it can only invoke a finite number of presuppositions. This fact is also a consequence of Gödel's Compactness theorem, which states that if a formula $\phi$ is a logical consequence of a possibly infinite collection of sentences $\mathcal{L}$, then it is a logical consequence of a finite subcollection of sentences of $\mathcal{L}$. The Compactness theorem can be derived from the Completeness theorem.} What does the incompleteness theorem say about this process? Briefly, that there cannot be a sound recursive collection of presuppositions from which to derive all mathematical truths. Any such collection will entail statements that are true in the intended domain but unprovable and others that may be assumed to be true or false safely and appended to the collection. An example of the former is $\sigma$, which is true in $\mathbb{N}$.\footnote{There are sentences that are independent of PA, are free of self-reference, and are more bona fide mathematical rather than metamathematical statements. The first such examples are perhaps with Paris-Harrington and by Kirby-Paris. Many excellent references on concrete incompleteness are available - for a recent one see Cheng 2021.} Examples of the latter include the Continuum Hypothesis, which is independent of Zermelo-Fraenkel set theory with the Axiom of Choice, and the Axioms of Choice and Determinacy that are contradictory to each other, and either of which but not both can be appended to  Zermelo-Fraenkel set theory. To some mathematicians this state of affairs seems quite natural. To others, notably Hilbert, it was a shocking realization. For clarity, it may be instructive to ask what presuppositions were required to prove the incompleteness theorem. A key one was that PA $\not \vdash S0 = 0$:  by assuming either PA $\vdash \sigma$ or PA $\vdash \lnot \sigma$ (in the Rosser version), $S0=0$ can be shown if the domain is fixed to be $\mathbb{N}$. The entire argument can be codified in PA. In light of incompleteness, first-order statements about numbers or other objects that can be codified in numbers have objective truth values but we will never have a theory that decides all such statements. Beyond that, given the independence of obviously cogent and intuitive propositions such as the Continuum Hypothesis, mathematical truth of infinite objects is subjective whilst the structure of mathematical argumentation is objective because it maps to first-order arithmetic and to computation.\footnote{I thank Scott Aaronson for pointing to me the distinction between objective truth of arithmetic statements and set-theoretic statements.}

\textbf{Computation.} Both the First and Second Incompleteness theorems can best be appreciated when the specific axiomatic system (PA, set theory), scope/power ($\omega$-consistency / soundness / simple consistency), proof strategy (diagonalization, Kolmogorov complexity, Turing computability) and proof implementation are all abstracted away. They are both algorithms that map across specific languages and domains, and which can be implemented in a multitude of ways. The proofs are literally computer code. In my opinion, the most important philosophical implications of the theorems are in computation. The first incompleteness proof gives birth to the notion of universal computation: total recursive functions are represented - \textit{implemented} - in any sufficiently rich axiomatic system, which is analogous to a programming language. Every computation is a deduction and every deduction is a computation - see also Kripke on the Church-Turing Thesis (2013). I note, however, that even though Turing machines, recursive functions, and deductions are mathematically equivalent in exhibiting universal computation, the notion misses practical aspects of important kinds of computation. Deep learning as well as neural computation in the human brain are impractical to map to recursive functions and what's worse, such mapping loses the inherent structure, symmetries, robustness, ability to train and dynamic nature that computational systems with high level of connection and analog gates exhibit. This limits the implications of the incompleteness results to philosophy of mind and AI, which brings us to the next topic.

\textbf{Human mind and AI.} Do the incompleteness results imply that human thought is superior to computation? The most famous argument is by Lucas-Penrose: any computational system is incomplete due to being precisely defined by some calculus; a human can prove that incompleteness; hence human reasoning cannot be mechanized. Many responses have been given by philosophers, mathematicians and computer scientists. In my opinion, the incompleteness results have no bearing on the ability of machines to match human thought. Part of the work of a mathematician is to formulate a collection of presuppositions relevant to the mathematical question in hand. Lets call this Process A, which is distinct from Process B of providing the argument. For much of ordinary mathematics,  Process A is reduced to picking an established domain, perhaps an open problem. Then Process B provides an argument that leverages a finite subset of presuppositions from the established axioms and theorems of that domain.\footnote{I am oversimplifying. The boundaries between Process A and Process B are not clear, and Process B usually involves creative definitions of new concepts, lemmas and corollaries. Those are important in the practice of mathematics, even though in principle Process B can be reduced to listing all syntactically valid proofs lexicographically: any theorem will eventually appear. In practice this is totally unreasonable, just like it is unreasonable to reduce a big neural network to a primitive recursive function or to a formula of PA.} Incidentally, the incompleteness theorem involves a nontrivial dose of Process A: (1) arithmetization is a trick to turn part of Process A into Process B: specifically, turn a deduction of $\phi$ from axiom instances $\phi_1, ..., \phi_k$ into a deduction of the deducibility of $\phi$; (2) in the final steps of the proof, PA $\not \vdash S0 = 0$ is assumed in order to prove that neither $\sigma$ nor $\lnot \sigma$ are provable; PA $\not \vdash S0 = 0$ is not part of PA. Back to the Lucas-Penrose argument: it assumes a computational system expressed purely as precise Process B calculus. However, Process A can be fallible and is not formalized in a precise calculus. Can it be mapped to a machine? I see no reason why not. For a language $\mathcal{L}$, $\sigma_{\mathcal{L}}$ can be deduced if the algorithm is allowed to assume $\mathcal{L} \not \vdash S0 = 0$. Whereas $\mathcal{L} \vdash \lnot \phi_{PROVABLE}(S0 = 0)$ leads to contradiction, $\mathcal{L} \not \vdash S0 = 0$ is an axiom that can be appended to $\mathcal{L}$ that never turns a sound $\mathcal{L}$ into an unsound one.\footnote{$\mathcal{L} + \{\mathcal{L}\not \vdash S0=0 \}$ is strictly stronger than $\mathcal{L}$, and appending "$\mathcal{L} + \{ \mathcal{L} \not \vdash S0=0 \} \not \vdash S0=0$" makes a stronger system, and so on. This process can continue. We can define ordinal $\omega$-many such extensions to create $\mathcal{L}_{\omega}$ as the union of all these axioms, and we can keep going even further to  $\mathcal{L}_{\omega} +\{\mathcal{L}_{\omega} \not \vdash S0=0 \}$ and so on, extending to Cantor's ordinals beyond infinity and producing ever stronger axiomatic systems. This is a beautiful justification of ordinals $>\omega$.} Programming Process A capability is challenging, but deep learning and other techniques will likely make progress in this direction. The incompleteness theorem says nothing against that. True, in principle a deductive system can express deep learning and other algorithms that may be deployed for Process A, but it can also express neural processes in the brain or even quantum system simulations that are sufficiently discretized and which deploy pseudorandom generators. Such mapping would be excruciatingly inefficient and physically impossible, and would require a ridiculous number of axioms to form $\phi_{PROVABLE}$. What makes minds work has to account for practical limitations. The Lucas-Penrose and similar arguments fail to do that, and the incompleteness results provide no support to the notion that human minds are superior to computation.

\section{Acknowledgements}
I thank Vassilis Gregoriades and Nick Nassuphis for insightful suggestions, edits and references, Scott Aaronson for comments that helped me clarify implications to mathematics, and Christos Athanasiadis for feedback.

\newpage

\section*{References}
\begin{enumerate}

\item Aaronson S. Shtetl-Optimized blog. Rosser's Theorem via Turing Machines. \\https://scottaaronson.blog/?p=710
\item Arjona M, Alonso E. Completeness: from Gödel to Henkin. History and Philosophy of Logic DOI: 10.1080/01445340.2013.816555, 2014.
\item Buldt B. The Scope of Gödel’s First Incompleteness Theorem. Logica Universalis (8):499-552, 2014.
\item Cheng Y. Current research on Gödel  incompleteness theorems. Bulletin of Symbolic Logic 27(2): 113-167, 2021.
\item Dean W. Recursive Functions. Stanford Encyclopedia of Philosophy, \\https://plato.stanford.edu/entries/recursive-functions/, 2020.
\item G\"{o}del, K. On formally undecidable propositions of Principia Mathematica and related systems I. Translated by Jean van Heijenoort. In Feferman, S. et al. (eds). Kurt G\"{o}del: Collected Works, Volume I (pp. 145–195). New York: Oxford University Press, 1986.
\item Hofstadter DT. G\"{o}del, Escher, Bach: An Eternal Golden Braid. Basic Books, 1979.
\item Kim B. Complete proofs of G\"{o}del's incompleteness theorems. https://web.yonsei.ac.kr/bkim/goedel.pdf.
\item Kikuchi M. Kolmogorov complexity and the second incompleteness theorem. Arch Math Logic 36, 437-443, 1997.
\item Krichman S, Raz R. The surprise examination paradox and the second incompleteness theorem. Notices of the AMS 57 (11), p. 1454-1458, 2010.
\item Kripke SA. The Church-Turing "Thesis" as a Special Case of Gödel's Completeness Theorem. In Computability: Turing, Gödel, Church and Beyond, BJ Copeland, C Posy and O Shagrir (eds). The MIT Press (Cambridge) 2013.
\item Kripke SA. The Road to Gödel. In: Berg J. (eds) Naming, Necessity, and More. Palgrave Macmillan, London.  2014.
\item Kripke SA. G\"{o}del's theorem and direct self reference. arXiv:2010.11979, 2021.
\item Nagel E, Newman JR. G\"{o}del's proof. NYU Press; Revised ed. 2001.
\item Smith P.  An Introduction to Gödel's Theorems, 2nd ed, Cambridge UK,  2013.
\item Smullyan R. Gödel's Incompleteness Theorems, Oxford Univ.Press, 1991.
\item Swierczkowski, S. Finite sets and Gödel’s incompleteness theorems. Dissertationes Mathematicae, 422, 1–58. 2003.
\item Weaver N. Forcing for mathematicians. World Scientific Publishing Co. 2014.
\end{enumerate}

\end{document}